\documentclass[11pt,english]{smfart}
\usepackage[cmtip,arrow]{xy}
\usepackage{pb-diagram,pb-xy}
\usepackage[dvipdfm]{graphicx}

\usepackage{epsfig}
\usepackage{amsmath}
\usepackage{amssymb}
\usepackage{amscd}
\usepackage{latexsym}
\usepackage{tabularx}
\usepackage{a4wide}
\usepackage{color}
\usepackage{enumerate}
\usepackage{subfigure}

\usepackage{url}
\usepackage{graphics}
\usepackage{mathrsfs}

\newtheorem{theorem}{Theorem}
\newtheorem{corollary}[theorem]{Corollary}
\newtheorem{proposition}[theorem]{Proposition}
\newtheorem{lemma}[theorem]{Lemma}

\newtheorem{definition}[theorem]{Definition}
\newtheorem{remark}[theorem]{Remark}
\newcommand{\dist}{\mathrm{dist}}
\newcommand{\ord}{\mathrm{ord}}
\newcommand{\D}{\mathrm{Div}}
\newcommand{\Red}{\mathrm{Red}}

\renewcommand{\int}{\mathrm{int}}
\renewcommand{\div}{\mathrm{div}}
\newcommand{\Wr}{\mathrm{Wr}}
\renewcommand{\O}{\mathcal O}
\renewcommand{\S}{\mathfrak S}

\title{Reduced divisors and Embeddings of Tropical Curves}
\author{Omid Amini}
\address{CNRS-DMA, \'Ecole Normale Sup\'erieure, 45 Rue d'Ulm, 75230 Paris Cedex 05}

\begin{document}
\begin{abstract}
Given a divisor $D$ on a tropical curve $\Gamma$, we show that reduced divisors define an integral affine map from the tropical curve to the complete linear system $|D|$. This is done by providing an explicit description of the behavior of reduced divisors under infinitesimal modifications  of the base point. We consider the cases where the reduced-divisor map defines an embedding of the curve into the linear system, and in this way,  classify all the tropical curves with a very ample canonical divisor. 
As an application of the reduced-divisor map, we show the existence of Weierstrass points on tropical curves of genus at least two and present a simpler proof of a theorem of Luo on rank-determining sets of points. We also discuss the classical analogue of the (tropical) reduced-divisor map: For a smooth projective curve $C$ and a divisor $D$ of non-negative rank on $C$, reduced divisors equivalent to $D$ define a morphism from $C$ to the complete linear system $|D|$, which is described in terms of Wronskians. 
\end{abstract}

\maketitle

\section{Introduction}

Reduced divisors are one of the main tools in the study of linear systems on tropical curves. Our aim in this paper is to present results on the behavior of reduced divisors under a modification of the base point and on their locus in the underlying complete linear system. 
 After a general introduction on tropical curves and linear systems, providing  the basic definitions and results, we show in Section~\ref{sec:reducedmap} that reduced divisors linearly equivalent to a divisor $D$ define an integral affine map from $\Gamma$ to the complete linear system $|D|$. This is done by providing an explicit description of the behavior of  reduced divisors under infinitesimal changes of the base point. We use  this description in Section~\ref{sec:canonicalembeddings} to characterize very ample divisors as the ones defining en embedding of $\Gamma$ in $|D|$, and classify all the tropical curves with a very ample canonical divisor.  In Section~\ref{sec:weier}, we present two applications of the reduced divisor map. We first show the existence of tropical Weierstrass points on tropical curves of genus at least two, and  second, we derive a simple proof of a theorem of Luo  on rank-determing sets of points on tropical curves. 
 In Section~\ref{sec:maps}, we compare the reduced divisor map with maps to tropical projective spaces defined by families of rational functions. Finally, in Section~\ref{app:weier} we show that the reduced divisor map has a classical analogue. For a smooth projective curve $C$ and a divisor $D$ of non-negative rank on $C$, reduced divisors equivalent to $D$ define a morphism from $C$ to the complete linear system $|D|$. This map is described in terms of Wronskians.

\subsection{Tropical curves and Riemann-Roch theorem}
 Tropical curves are algebraic curves defined over the tropical semi-ring $\mathbb T = (\mathbb R\cup \{-\infty\},\oplus,\odot)$. The tropical sum $\oplus$ is taking the maximum and the tropical product $\odot$ is the (usual) sum, i.e., $a\oplus b =\max\{a,b\}$ and $a\odot b = a+b$. 
 A (complete) tropical curve is a compact topological space homeomorphic to a one-dimensional simplicial complex equipped with an integral affine structure. This means that the topological space has a finite covering by open sets $U_i$ and there are chart maps $\phi_i: U_i \rightarrow \mathbb T^{n_i}$, sending $U_i$ homeomorphically to a subset $V_i \subset \mathbb T^{n_i}$ such that the change of charts are restrictions of affine maps with integer linear part, i.e., for $i$ and $j$, there exist a matrix $A^{i,j}\in M_{n_j\times n_i}(\mathbb Z)$ and an $a \in \mathbb T^{n_j}$ such that the change of charts $\phi_j\phi_i^{-1}$ on $\phi_i(U_i \cap U_j)$ coincide with the map $A.x+a$. 
It is easy to see that the above definition is equivalent to the following more concrete definition.

\vspace{.3cm }

A weighted graph $(G,\ell)$ is a finite graph $G = (V,E)$ where each edge $e\in E$ has been assigned a positive length $\ell_e \in \mathbb T$. This data defines a metric space $\Gamma$, called the geometric representation of $(G,\ell)$, as follows. This is obtained by replacing each edge of the graph by an  interval of the given length associated to that edge, the end points of the interval being the corresponding vertices of the original weighted graph. The metric space $\Gamma$ is called a metric graph and the corresponding weighted graph $(G,\ell)$ is called a model of $\Gamma$. It is clear that different models can give the same metric graph. The choice of a model corresponds to the choice of a finite set of vertices $V$, i.e., a subset of $\Gamma$ such that $\Gamma \setminus V$ is a disjoint union of open intervals. In this paper, when we talk about the vertices of a metric graph, we understand that a model of $\Gamma$ is fixed.

 Note that a metric graph $\Gamma$ inherits a topology from the above described geometric representation and it is possible to speak about the closed and open sets.  The distance $\dist_\Gamma (v, w)$ between two points $v$ and $w$ of $\Gamma$ is measured in the metric space corresponding to the geometric representation of $\Gamma$ (the subscript $\Gamma$ is omitted if the metric graph $\Gamma$ is clear from the context). A segment in $\Gamma$ is a closed subinterval of one of the intervals in a geometric representation of $\Gamma$ described above.
 
 Each point $v$ of $\Gamma$ has a degree $\deg(v)$. This is defined by taking a sufficiently small neighborhood of $v$ in $\Gamma$, removing $v$ from that neighborhood and counting the number of connected components. Note that every point in the interior of an edge of a model of $\Gamma$ has degree two. Points of degree larger than two are called branching points. The coarsest model for $\Gamma$, if it exists, is a weighted graph in which every vertex has degree  different from two, i.e., each vertex is either a branching point or a leaf.  Note that the coarsest model is unique if it exists. And it exists if the metric graph is connected and there is at least one branching point, i.e., if the metric graph is not a circle (a weighted graph is allowed to contain loops). 
 
 An abstract tropical curve is a metric graph where some of the edges incident with vertices 
of degree one (leaves) have infinite length. Such edges are identified with the 
interval $[0, \infty]$, such that $\infty$ is identified with the vertex of degree one, and 
are called infinite ends or bounds. Since infinite ends do not play any interesting role with respect to divisor theory on a tropical curve, in this paper we suppose that $\Gamma$ does not have any infinite end and so, in particular, our tropical curves are all compact. 

The genus of a tropical curve $\Gamma$, denoted by $g(\Gamma)$ or $g$, is the first Betti number of $\Gamma$ as a topological space. Given a model $(G,\ell)$ of $\Gamma$ with $G = (V,E)$, the genus of $\Gamma$ is  $g=|E|-|V|+1$.

 \vspace{.3cm}

A divisor $D$ on a tropical curve $\Gamma$ is an element of the infinite free abelian group $\D(\Gamma)$ on points of $\Gamma$. The basis element of $\D(\Gamma)$ corresponding to a point $v$ of $\Gamma$ is denoted by $(v)$.  Each divisor $D \in \D(\Gamma)$ can be uniquely written as 
$D = \sum_{v\in\Gamma} a_v(v)$ with $a_v \in \mathbb Z$ and all but a finite number of $a_v$'s are zero. The coefficient $a_v$ of $(v)$ in $D$ is denoted by $D(v)$. A divisor $D$ is called effective,  if $D \geq 0$, i.e., if  $D(v) \geq 0$ for 
all $v \in \Gamma$.   The degree of a divisor $D \in \D(\Gamma)$, denoted by $\deg(D)$, is defined as 
$\deg(D) := \sum_{v\in\Gamma}D(v).$ 
The canonical divisor of a tropical curve $\Gamma$ is the divisor $K_\Gamma$ defined by
\[K_\Gamma :=\sum_{v\in \Gamma} (\deg(v)-2)(v).\]

A rational function on a tropical curve $\Gamma$ is a continuous piecewise linear function $f : \Gamma \rightarrow 
\mathbb{R} \cup \{\pm\infty\}$ with integral slopes at every point and such that the number of linear parts of $f$
on every edge is finite and the only points $v$ with $f (v) = \pm\infty$ are unbounded  ends. (This is simply equivalent to a continuous map respecting the integral affine structure of $\Gamma$ and the obvious integral affine structure on $\mathbb{R} \cup \{\pm\infty\}$.)  The order $\ord_f (v)$ of a rational function $f$ on $\Gamma$ at a point $v$  
is the sum of the outgoing slopes of $f$ along the segments of $\Gamma$ emanating from $v$. 
Given a rational function $f$, the  divisor $\div(f) = \sum_v \ord_v(f)(v)$ is called the principal divisor associated to $f$. Note that if $v$ is not a branching point of $\Gamma$ and the function $f$ does 
not change its slope at $v$, then $\ord_f (v) = 0$, thus, there are only finitely many 
points $v$ with $\ord_f (v)\neq 0$ and $\div(f) \in \D(\Gamma)$. The degree of a principal divisor is equal to zero since each linear part of given slope in a rational function contributes twice in the degree, once with positive sign and once with negative sign.  Two divisors $D$ and $D'$ on 
$\Gamma$ are called linearly equivalent, written $D \sim D'$, if there exists a rational function $f$ on $\Gamma$  such that $D = D'+ \div(f)$. 

Let $D$ be a divisor on a tropical curve $\Gamma$. The (complete) linear system of $D$, denoted by $|D|$, is defined as the set of all effective divisors $E$ linearly equivalent to $D$, 
\[|D|:= \{\,E\in\D(\Gamma)\::\: E\geq 0 \textrm{ and }E\sim D\,\}.\]
The rank of a divisor $D \in \D(\Gamma)$, denoted by $r(D)$, is defined by
 \[r(D):= \min_{E\::\:E\geq 0, |D-E|=\emptyset}\deg(E)-1.\]
 
 \vspace{.3cm}

Gathmann and Kerber~\cite{GK08} and, independently, Mikhalkin and Zharkov~\cite{MZ08}
proved the following Riemann-Roch theorem for tropical
curves (see also~\cite{HKN08} for another more combinatorial proof of this theorem similar to the proof of Baker and Norine's Riemann-Roch theorem for graphs~\cite{BN06}).

\begin{theorem}[Tropical Riemann-Roch]\label{thm:R-R} Let $D$ be a divisor on a tropical curve $\Gamma$ of genus $g$ and $K$ be the canonical divisor of $\Gamma$. Then
\[r(D) - r(K- D) = \deg(D) + 1 - g.\]
\end{theorem}

\vspace{.5cm}

\subsection{Reduced divisors}
The proof of Riemann-Roch theorem for graphs and tropical curves is based on the notion of reduced divisors that we now explain. (We refer to~\cite{AM09} for a different, more geometric proof, of a slightly more general result in the graph case.) 

\vspace{.3cm}

Let $\Gamma$ be a tropical curve and $v_0$ be a given (base) point of $\Gamma$. 
For a closed connected subset $X$ of $\Gamma$ and a point $v \in \partial X$, the number of edges leaving $X$ at $v$, denoted by $\deg^{out}_X(v)$, is by definition the maximum size of a
collection of internally disjoint segments in $\Gamma \setminus (X -\{v\})$ with an end in $v$. In what follows we also call such a connected subset of $\Gamma$ a {\it cut} in $\Gamma$.
A boundary point $v$ of a closed connected subset $X$ of  $\Gamma$ is called
{\it saturated} with respect to $D \in \D(\Gamma)$ and the cut $X$ if the number of edges
leaving $X$ at $v$ is at most $D(v)$ (i.e., if $\deg^{out}_X(v) \leq D(v)$) and is called {\it non-saturated}, otherwise. 
A cut is said to be saturated if all its boundary points are saturated.  
When talking
about saturated and non-saturated points, we omit sometimes $D$ and/or the set $X$ if they are clear from the context. 
 \vspace{.3cm}
 
 A divisor $D$ on a tropical curve $\Gamma$
is said to be {\it $v_0$-reduced} if the following two properties are satisfied: 
\begin{itemize}
\item For all  point $v \neq v_0$  of $\Gamma$, $D(v) \geq 0$, in other words if all the coefficients of $D$ are non-negative except possibly for the base-point $v_0$. 
\item For every closed connected subset $X$ of points of  $\Gamma$ which does not contain $v_0 $, there exists a non-saturated point $v \in \partial X$.
\end{itemize}

\vspace{.3cm}

 Every divisor is equivalent to a unique $v_0$-reduced divisor. We provide here a quick proof and refer to~\cite{BN06, HKN08} for more details.
\begin{theorem}
Let $\Gamma$ be a tropical curve and $v_0$ be a point of  $\Gamma$. For every divisor $D \in \D(\Gamma)$, there exists a unique $v_0$-reduced divisor $D_0$ such that $D_0 \sim D$.
\end{theorem}
\begin{proof}

Let $\mathcal T_D$ be the set of all divisors $D'\sim D$ such that all the coefficients of $D'$ are non-negative at every point of $\Gamma$ except possibly at $v_0$, and such that, in addition, the coefficient of $v_0$ is maximum possible with respect to this property. Note that $\mathcal T_D$ is non-empty. This is simply because there are linearly equivalent divisors to $D$  with non-negative coefficients at every point $v \neq v_0$ of $\Gamma$. To see this, note that if $k$ is large enough (as a function of the number of branching points, the diameter of the underlying graph, and the maximum absolute value of the coefficients of $D$), then it can be easily shown that the divisor $D+k (v_0)$ is equivalent to an effective divisor $E$. Now $E -k (v_0)\sim D $ has the required property. Note that Riemann-Roch theorem implies that if $k$ is large enough so that the degree of $D+k(v_0)$ is larger than $g$, then the rank of $D+k(v_0)$ is non-negative, however the existence of a large enough $k$ with this property does not need the full strength of the Riemann-Roch theorem (whose proof is bases on reduced divisors).

\vspace{.3cm}

We next observe that $\mathcal T_D$ inherits a natural topology from the topology of $\Gamma$ and is compact.

Let $N$ be the sum of the coefficients of points of $\Gamma\setminus\{v_0\}$ in any divisor in $\mathcal T_D$ (by the definition of $\mathcal T_D$, $N$ is well-defined). Let $\mathcal A$ be the subset of $\mathbb R^{N}$ defined as follows
$$\mathcal A := \Bigl\{(x_1,\dots,x_N)\,|\,0\leq x_1\leq x_2\leq \dots\leq x_N\Bigr\}.$$
Note that $\mathcal A$ has a natural total order defined by the lexicographical order: $(x_1,\dots,x_N) < (y_1,\dots,y_N)$ if and only if there exists a non-negative integer $i$ such that $x_1=y_1,\dots,x_i =y_i$ and $x_{i+1} < y_{i+1}$.

Define a continuous map $F: \mathcal T_D \rightarrow \mathcal A$ as follows: for each divisor $D' \in \mathcal T_D$, consider the multiset $A(D')$ of points in $\Gamma \setminus \{v_0\}$ where each point $v\neq v_0$ appears in this multiset exactly $D'(v)$ times. Define $F(D')$ to be the point of $\mathcal A$ defined by the multiset of distances $\dist_\Gamma(v,v_0)$ for $v\in A(D')$ ordered in an increasing way to define a point of $\mathcal A$. It is straightforward to show that the map $F$ is continuous.

Since $\mathcal T_D$ is a compact set and $F$ is  continuous,  there exists a divisor $D_0$ such that $F$ takes its minimum value on $D_0$, i.e., $F(D_0) = \min_{D'\in \mathcal T_D}F(D')$. We claim that $D_0$ is $v_0$-reduced. The first property is clearly verified by the choice of $D_0$: $D_0 \sim D$ and the coefficients $D_{0}(v)$ are all non-negative for $v \neq v_0$ in $\Gamma$. We must show that every connected set $X$ which  does not contain $v_0$ has a non-saturated point on its boundary. For the sake of a contradiction, suppose this is not the case and let $X$ be the connected set violating this condition. This means that for all $v \in \partial X$, $D_{0}(v) \geq \deg^{out}_X(v)$. By the definition of $\deg^{out}_X$, there exists an $\epsilon>0$ such that for each vertex $v \in \partial X$, there exist closed segments $I_1^v,\dots,I^v_{\deg^{out}_X(v)}$ emanating from $v$ with the following properties:
\begin{itemize}
\item For $v \in \partial X$ and $1\leq  j\leq \deg^{out}_X(v)$, each $I_j^v$ has length $\epsilon$,  and the half-open line segments $I^v_j \setminus v$ live outside $X$ and are all disjoint.
\end{itemize}
Define now the rational function $f: \Gamma \rightarrow \mathbb R$ as follows: $f$ takes value zero on $X$, is linear of slope $-1$ on each interval $I_j^v$ and takes value $-\epsilon$ on all the points of $\Gamma \setminus \Bigl(X\cup \bigcup I_j^v\Bigr)$, the union being over all $v\in \partial X$ and $1\leq j\leq \deg^{out}_X(v)$. Since $\deg^{out}_X(v) \leq D_{0}(v)$, the divisor $D_0 + \div(f)$ lies inside $\mathcal T_D$. Since the distance of $X$ to $v_0$ is strictly larger than the distance from the set $X \cup \bigcup I^v_j$ to $v_0$, and the boundary points of $X$ are all in the multiset $A(D_0)$, it is straightforward to see that $F(D_0) > F(D_0 + \div(f))$, and this contradicts the choice of $D_0$. This proves the existence. 

To prove the uniqueness, for the sake of a contradiction, suppose there are two $v_0$-reduced divisors $D_0$ and $D_0'$ which are linearly equivalent. This means that there exists a non-constant rational function $f$ such that $D'_0 = D_0+\div(f)$. Since $f$ is not a constant function, we can assume that $f$ does not take its maximum at $v_0$ (if not, we change the role of $D_0$ and $D'_0$, and replace $f$ by $-f$). Let $X$ be a connected component of the set of all points where $f$ takes its maximum. Note that $v_0 \notin X$. For any point $v$ lying on the boundary of $X$, the slope of $f$ along any segment $I_j^v$ emanating from $v$ is at most $-1$. It follows that the coefficient of $(v)$ in $D_0 + \div(f)$ is at most $D_0(v) -\deg^{out}_X(v)$. Since $D_0$ is $v_0$-reduced, there exists a point $v\in \partial X$ such that $D(v) < \deg^{out}_X(v)$. For this point, $D'_0(v) = D_0(v) - \deg^{out}_X(v) <0$, which contradicts the assumption that $D'_0$ is non-negative outside  $v_0$. The uniqueness follows.
\end{proof}

\section{The Map $\Red$ From $\Gamma$ to a (Complete) Linear System}\label{sec:reducedmap}

Let $D$ be a divisor and $|D|$ its complete linear system. Suppose that $r(D)\geq 0$ such that $|D|$ is non-empty. There is a natural polyhedral structure on $|D|$ that we now explain.

Let $\Gamma$ be a tropical curve and $N$ be an integer. The symmetric product $\Gamma^{(N)}$ is defined as the set of all unordered $N$-tuples of points of $\Gamma$, i.e.,
\[\Gamma^{(N)}:= \Gamma^N / \S_N,\]
where $\S_N$ is the symmetric group of degree $N$ and $\Gamma^N$ is the product of $N$ copies of $\Gamma$.

\noindent The linear system $|D|$ is naturally a subset of the symmetric product  $\Gamma^{(\deg(D))}$:
define a map $\iota: |D| \hookrightarrow \Gamma^{(\deg(D))}$ by sending each $E \in |D|$ to $\iota(E)$ of $\Gamma^{(\deg(D))}$ consisting of the (unordered) multiset of all the points $v$ with $E_v\neq 0$.
The topology of $|D|$ is the natural topology induced from $\Gamma^{(\deg(D))}$. We will explain  the cell structure of $|D|$ in more details in the next section.

\subsection{The integral affine structure on $\Gamma^{(N)}$ and $|D|$} \label{sec:integeraffine} We start this section by describing a natural integral affine structure on $\Gamma^{(\deg(D))}$  induced from $\Gamma$. First, we need to say more precisely what we mean by an integral affine structure. Recall that a manifold with an integral affine structure of dimension $n$ is given by the following data:
\begin{itemize}
\item A (finite) collection of open sets (charts) $U_i$ covering all of $X$. 
\item A homeomorphism $\phi_i$ from $U_i$ to an open subset $V_i$ of $\mathbb R^{n}$ for each $i$ such that for different $i$ and $j$ the change of chart map $\phi_j\,\phi_i^{-1}: \phi_i( U_j\cap U_i) \rightarrow \phi_j(U_i\cap  U_j)$ is the restriction of an integral affine map from $\mathbb R^n$ to $\mathbb R^n$, i.e., there exists $a \in \mathbb R^{n_j}$ and $A \in M_{n\times n}(\mathbb Z)$ such that $\phi_i\,\phi_j^{-1}(x) = Ax+a$ for all $x \in \phi_i(U_j\cap U_i)$.
\end{itemize} 
Since the objects we will be considering will have singularities, we need to modify  the above definition to cover also the case of stratified manifolds. Roughly speaking,  if $X$ admits a filtration $X_0 \subset X_1\subset \dots \subset X_{n-1} \subset X_n = X$ such that for each $i$, $X_{i} \setminus X_{i-1}$  is a (possibly disconnected) manifold of dimension $i$, then it is natural to impose that each $X_i\setminus X_{i-1}$ is an $i$-dimensional manifold with an integral affine structure and such that in addition, there is a compatibility between the integral affine structure on $X_i\setminus X_{i-1}$ and the one "induced" from $X_{i+1}\setminus X_i$ on $X_i\setminus X_{i-1}$.  To give a more precise meaning to this, in what follows it turns out to be more convenient  for our purpose to work with abstract rational polyhedral complexes (which will carry by definition an integral affine structure in the above general sense).

Let $N$ be a lattice of rank $n$ in the $n$-dimensional real vector space $N_\mathbb R = N \otimes_\mathbb Z \mathbb R$. Recall that a rational polyhedral complex $\Sigma$ in $N_\mathbb R$ is a finite collection of (convex) rational polyhedra in $N_\mathbb R$ which verifies the following two properties:
\begin{itemize}
\item[$\bullet$] If $\sigma \in \Sigma$ and $\tau$ is a face of $\sigma$, then $\tau$ is in $\Sigma$;

\item[$\bullet$] If $\sigma$ and $\tau$ are in $\Sigma$, then $\sigma\cap \tau$ is a (possibly empty) common face of $\sigma$ and $\tau$ (which is also in $\Sigma$).  
\end{itemize}

Since the polyhedral complex $\Sigma$ is rational, for each point $x$ in the relative interior of a polyhedra $\sigma \in \Sigma$, the tangent space of $x$ in $\sigma$ has a canonical sublattice of full rank. This  canonically coincides with the sub-lattce $N_\sigma$ of $N$ of dimension $\dim(\sigma)$ defined by the set of all points of $N$ which lie on the rational plane $\mathrm{span}(\sigma - x)=:N_{\sigma,\mathbb R}$ in $N_\mathbb R$ defined by $\sigma$.  Under this identification, one canonically has an isomorphism between the tangent space of the interior of $\sigma$ with $\int(\sigma)\times N_{\sigma,\mathbb R}$.  

Now to make the above definition independent of the embedding into a fixed $N_\mathbb R$, we can proceed as follows. An abstract rational polyhedral complex structure $\Sigma$ on a Hausdorff topological space $X$ is a finite collection of closed topological subspaces $\sigma \subseteq X$ with the following properties. Each $\sigma$ is homeomorphic to a full dimensional polyhedra $|\sigma|$ in a real vector space $N_{\sigma,\mathbb R}$. In addition, $N_{\sigma,\mathbb R}$ has a sublattice of full rank $N_\sigma$ such that $|\sigma|$ is rational with respect to $N_\sigma$. By rational, we mean that the real vector space defined by any face of $|\sigma|$ comes (by extension of scalars from $\mathbb Q$ to $\mathbb R$) from a sub-vector space of $N_{\sigma,\mathbb Q} = N\otimes_\mathbb Z \mathbb Q$.  Under these homeomorphisms, $\Sigma$ verifies the following properties, analogue to the properties of a rational polyhedral complex in $N_\mathbb R$ given above. For $\sigma \in \Sigma$, define $F_\sigma$ as the poset (under inclusion) of all $\tau \in \Sigma$ such that $\tau \subseteq \sigma$. For convenience, assume that $\emptyset \in F_\sigma$ for all $\sigma \in \Sigma$.

\begin{itemize}
\item[$\bullet$] If $\sigma \in \Sigma$,  there is an isomorphism of posets from $F_\sigma$ to the face poset of $|\sigma|$.  Under this isomorphism, if $\tau \in F_\sigma$ and $\eta$ is the corresponding face of $|\sigma|$ in $N_{\sigma,\mathbb R}$,  there exists an isomorphism of rational polyhedra between $|\tau|$ and $\eta$. 
 
\item[$\bullet$] If $\sigma$ and $\tau$ are in $\Sigma$, then $\sigma\cap \tau$ is in $\Sigma$ and thus belongs to both $F_\sigma$ and $F_\tau$.
\end{itemize}

Given an abstract rational polyhedral complex structure $\Sigma$ on $X$, define $X_i$ as the union of all $\sigma \in \Sigma$ such that $|\sigma|$ has dimension $i$. The following facts are easy to verify:
\begin{itemize}
\item For each $i$, $X_i \setminus X_{i-1}$ is naturally a (open) manifold of dimension $i$ with an integral affine structure. More precisely, $X_i \setminus X_{i-1}$ is the disjoint union of the relative interior of all $\sigma \in \Sigma$ such that $|\sigma|$ has dimension $i$.
 
 \item If $\Sigma$ and $\Sigma'$ are abstract rational polyhedral structures on $X$ and $X'$, respectively, there exists a natural abstract rational polyhedral complex structure $\Sigma\times \Sigma'$ on $X \times X'$. The elements of $\Sigma\times \Sigma'$ are of the form $\sigma\times \sigma' \subset X\times X'$ for $\sigma \in \Sigma$ and $\sigma'\in \Sigma'$.
\end{itemize}

Any metric graph $\Gamma$ with a simple model $G=(V,E)$ admits an abstract rational polyhedral structure given by $\Sigma = V \cup E$. (Recall that $G$ is called simple if it does not have any loop or parallel edges.)  For each edge $e \in E$, identify $e$ with the interval $[0,\ell(e)]$ of length $\ell(e)$ in $\mathbb R$ and let $N_e = \mathbb Z \subset \mathbb R$. This naturally defines an abstract rational polyhedral complex structure $\Sigma^N = \Sigma^N_{G=(V,E)}$ on $\Gamma^N$.  Each $\sigma \in \Sigma_{G=(V,E)}$ is of the form $\sigma_1\times\dots\times \sigma_N$ for $\sigma_i\in V\cup E$.  
Note that the integral affine structure on $\Gamma^N$ is invariant under the action of the symmetric group $\mathfrak S_N$. Note also that subdividing the rational polyhedral structure of $\Gamma$ consists in  choosing another model $G'=(V',E')$ of $\Gamma$ which is a refinement of $G$ (i.e., $V \subseteq V'$). This in turn induces a subdivision of the abstract rational polyhedral structure of $\Gamma^{N}$.

\medskip

We next show that the topological quotient $\Gamma^{(N)}$ of $\Gamma^N$ by $\S_N$ admits an integral affine structure given by an abstract rational polyhedral complex structure $\Sigma^{(N)}$. In addition, the projection map $\pi:\Gamma^N \rightarrow \Gamma^{(N)}$ is "integral affine " in the sense that for all $x$ in the relative interior of a $\sigma$ in $\Sigma^N$, with image $y$ in the relative interior of $\eta$ in $\Sigma^{(N)}$, $N_\eta \subseteq \pi_*(N_\sigma)$ ($\pi_*$ is the push forward of tangent vectors from $\Gamma^N$ to $\Gamma^{(N)}$).

To describe the integral affine structure of $\Gamma^{(N)}$, we fix a model $G=(V,E)$ of $\Gamma$ without loops and parallel edges. From the description below, it will be clear that choosing another simple model $G'=(V',E')$ of $\Gamma$ which is a refinement of $G$ consists in subdividing the abstract rational polyhedral structure of $\Gamma^{(N)}$, and thus, the induced integral affine structure does not depend on the particular choice of the model.

\noindent  Let $v_1,\dots,v_n$ be the set of all the vertices of $V$. The enumeration defines an orientation on the set of edges of $\Gamma$, thus a lexicographical order on the points of the edge respecting the orientation, and also a lexicographical order on the edges set. Since $\Gamma$ does not have any loop or parallel edges, there is no ambiguity in the definition of the lexicographical order.
First, for an edge $e$ of $G$ of length $\ell(e)$ and an integer $m$, consider the closed subset $e^{(m)} = e^m/\mathfrak S_m$ of $\Gamma^{(m)}$, and indentify $|e^{(m)}|$ with the rational polytope in $\mathbb R^m$ consisting of all the points $(x_1,\dots,x_m)$ such that $0\leq x_1\leq x_2\leq \dots\leq x_m \leq \ell(e)$. The faces of $|e^{(m)}|$ induce a decomposition of $e^{(m)}$ into closed subsets and define a rational polyhedral structure on $e^{(m)}$. 
Now, more generally, given a sequence $\underline e$ consisting of $N$ edges $e_1 =e_1\dots=e_1< e_2 = \dots =e_i < \dots < e_k = \dots= e_k$ of $\Gamma$ ordered in the lexicographical order, and in which each edge $e_i$ is repeated $s_i$ times with $\sum_{i=1}^k\, s_i =N$, consider the closed subset $\sigma_{\underline e}$ of $\Gamma^{(N)}$ defined by $\sigma_{\underline e} = e_1^{(s_1)} \times \dots \times e_k^{(s_k)}$, and consider its geometric realization as a rational polytope  in $\mathbb R^N$ defined by 
\[|\sigma_{\underline e}| = |e_1^{(s_1)}| \times \dots \times |e_k^{(s_k)}|.\]
The faces of $|\sigma_{\underline e}|$ define a decomposition of $\sigma_{\underline e}$ into closed subsets, and induce an abstract rational polyhedral structure $\Sigma_{\underline e}$ on $\sigma_{\underline e}$.  Note that the relative interior $\int(\sigma_{\underline e })$ of $\sigma_{\underline e}$ is given by
\begin{align*}
\int(\sigma_{\underline e})= \Bigl\{(x_1,\dots,x_{s_1}, x_{s_1+1},\dots,&x_{s_1+s_2},\dots, x_{s_1+\dots+s_{k-1}+1},\dots,x_N) \in \Gamma^N\,|
\\
&\, x_1< \dots < x_{s_1} \in \int(e_1)\,;\\
&\,x_{s_1+1} < \dots < x_{s_1+s_2} \in \int(e_2)\,;\,\dots;\\
&\,x_{s_1+\dots+s_{k-1}+1}< \dots < x_N \in \int(e_k)\Bigr\}.
\end{align*}

Here $\int(e)$ denotes the interior of an edge $e$.
The abstract rational polyhedral structure $\Sigma$ on $\Gamma^{(N)}$ is defined as follows.  $\Sigma^{(N)}$ is the union of all the closed subsets of $\Gamma^{(N)}$ which appear in some $\Sigma_{\underline e}$ for some sequence of edges $\underline e$, $\Sigma^{(N)} = \bigcup_{\underline e} \Sigma_{\underline e}$. Each element $\tau$ of $\Sigma^{(N)}$ belongs to some $\Sigma_{\underline e}$, and thus one can define $|\tau|$ and its rational polyhedral structure to be the one defined by the corresponding face of $|\sigma_{\underline e}|$. It is straightforward to check that this does not depend on the choice of the sequence $\underline e$, and thus one obtains a global structure of an abstract rational polyhedral complex on $\Gamma^{(N)}$.   This abstract rational polyhedral structure on $\Gamma^{(N)}$ is pure of dimension $N$, i.e., all the maximal cells in $\Sigma^{(N)}$ have dimension $N$.

Note that in the stratification $X_0 \subset X_1\subset \dots \subset X_N$ of $\Gamma^{(N)}$ defined by $\Sigma^{(N)}$, $X_N \setminus X_{N-1}$ is the disjoint union of all the open sets $\int(\sigma_{\underline e)}$. If among the strict inequalities describing $\int(\sigma_{\underline e})$ $i$ of them become equality (counting also the situations where the points lie on the boundary of an edge), then these so defined cells become disjoint and their union form the subset $X_{N-i} \setminus X_{N-i-1}$. It is clear that the set of vertices $X_0$ consists of all the $N$-tuples $v_{i_1},\dots,v_{i_N}$ of vertices of $G$ with $1\leq i_1\leq \dots \leq i_N \leq n$.

\vspace{.5cm}

Our next aim will be to explain how $|D|$ inherits a natural structure of an abstract rational polyhedral complex from that of $\Gamma^{(\deg(D))}$. This essentially amounts in providing the description of the cell structure on $|D|$ given in~\cite{GK08, MZ08} and~\cite{HMY09}. We point out that this cell structure is rational, is induced by the rational polyhedral complex structure of $\Gamma^{(\deg(D))}$, and is so that the inclusion map $|D| \hookrightarrow \Gamma^{(\deg(D))}$ is integral affine.

Let $\Gamma$ be a tropical curve with a simple model $G=(V,E)$. Let $v_1,\dots,v_n$ be the set of all the vertices of $V$. Let $D$ be a divisor of degree $d=\deg(D)$ of non-negative rank.
Consider the abstract rational polyhedral complex structure $\Sigma^{(d)}$ on $\Gamma^{(d)}$ described above. For each $\tau \in \Sigma^{(d)}$, the intersection of $|D|$ with the relative interior of $\tau$ has a finite number of connected components.  Define $\Sigma_D$ as the family  of all the closed subsets of $|D|$ which are the closure in $\Gamma^{(d)}$ (or equivalently in $|D|$) of a connected component of $|D| \cap \int(\tau)$ for a $\tau \in \Sigma^{(d)}$. We claim that $\Sigma_D$ defines an abstract rational polyhedral structure on $|D|$. 
This is reduced to proving two claims, namely that first, for each $\eta$ the closure of a connected component of $|D| \cap \int(\tau)$, 
the geometric realization $|\eta|$ in $|\tau|$ is a rational polytope, and that in addition, each face of $|\eta|$ is the geometric realization of an element $\gamma$ in $\Sigma_D$, this second assertion being itself a direct consequence of the first one and the assertion that the dimension of $|D| \cap \int(\tau)$ is locally constant (thus, $\eta$ cannot have a proper face in the interior of $\tau$).

 Recall that the enumeration of the vertices induce an orientation of the edges of $\Gamma$ and a lexicographical order on the edges. The relative interior of an element $\tau \in \Sigma^{(d)}$ is the set of all the points $\underline x = \{x_1,\dots,x_d\}$ of $\Gamma^{(d)}$ described by 
 \begin{itemize}
\item A map $\mathfrak d: V \sqcup E \rightarrow \mathbb Z$ which encodes the following data: For any vertex $v\in V$, $\mathfrak d(v)$ is the number of $x_i$ in $\underline x$ which are equal to $v$. For any edge $e\in E$, $\mathfrak d(e)$ is the number of $x_i$ in $\underline x$ which lie in the interior of the edge $e$. Note that $\sum_{v\in V} \mathfrak d(v) + \sum_{e\in E} \mathfrak d(e) =d$.
\item In addition, for every edge $e$, an ordered decomposition of $\frak d(e) = \frak d_1(e)+\dots+\frak d_{s_e}(e)$ for some non-negative integer $s_e=s_e(\tau)$. This is given by looking at the location of all the points $x_i$ which are in the interior of the edge $e$: among these points, the first $\mathfrak d_e(1)$ are the same and are strictly smaller than the following $\mathfrak d_2(e)$ points which are the same, these points being strictly smaller than the following $\mathfrak d_3(e)$ points which are the same, and so on.
\end{itemize}
 
 We now claim that $\int(\tau) \cap |D|$ is a disjoint union of rational polyhedra in $N_{\tau,\mathbb R}$ of the same dimension. This establishes the two above mentioned assertions and thus show that $|D|$ has an abstract rational polyhedral complex structure.  

\noindent To see this, first note that $|D|\cap \int(\tau)$ consists of  all the divisors $D'$ linearly equivalent to $D$ such that 
\begin{itemize}
\item the coefficient of $D'$ at $v$ is equal to $\mathfrak d(v)$, and
\item for every edge $e$, there are points $x_1^e<x_2^e<\dots< x_{s_e}^e \in \int(e)$ such that the support of $D'$ in the interior of $e$ is given by $\sum \mathfrak d_i(e)(x_i^e)$. 
\end{itemize} 
To describe the rational polyhedral structure of $|D| \cap  \int(\tau)$, one also has to fix the slopes of the rational function $f$ which provides the linear equivalence of $D$ and $D'$, i.e., $ \div(f)+D =D'$, at the beginning of any edge $e$. In other words, one fixes an integer valued function $m: E \rightarrow \mathbb Z$ and assumes that $f$ has slope equal to $m(e)$ along $e$ at the starting point of $e$. The combinatorial data $\mathfrak d(.)$, $\mathfrak d_i(.)$ and $m(.)$ provides linear equations and inequalities with rational coefficients which define a rational polyhedral part of $|D| \cap \int(\tau)$~\cite{GK08}. In addition, by the up-to-an-additive-constant uniqueness of the rational function $f$, it becomes immediately clear that these polyhedra form  different connected components of $|D| \cap \int(\tau)$. To show the assertion that the dimensions are all equal,   let $n_\tau$ be the number of connected components of the topological space  obtained by removing $s_e$ points from the interior of any edge $e$. Note that for any $D' = \sum_{v \in V}  \mathfrak d(v)(v)  + \sum_{e\in E} \sum_{i=1}^{s_e} \mathfrak d_i(e) (x^e_i)\in |D|\cap \int(\tau)$, the number of connected components of $\Gamma \setminus \{x^e_i\}_{i,e}$ is $n_\tau$. Each connected component $A$ of $\Gamma \setminus \{x^e_i\}_{i,e}$
 forms a saturated cut in $\Gamma$, because all the points on the boundary have out-going degree one and coefficient at least one. For small enough positive $\epsilon$, by considering small out-going segments $I_\epsilon^v$ of length $\epsilon/D'(v)$ from $A$ at $v$ for any point $v\in \partial A$, one can define a rational function $f_\epsilon$ with value zero at $A$, with slope $-D'(v)$ along the segment $I_\epsilon^v$, and with value $-\epsilon$ outside the union of $A$ and the segments $I_\epsilon^v$. In this way, by combining these functions, one can show that the $n_\tau$ connected components  of $\Gamma \setminus \{x^e_i\}_{i,e}$ allow to infinitesimally move $D'$  in an open subset of dimension $n_\tau -1$ in $|D| \cap \int(\tau)$, which shows that the dimension of each connected component of $|D| \cap \int(\tau)$ is $n_\tau-1$, c.f.~\cite{HMY09} and Proposition~\ref{prop:HMY}.

\medskip

 It can happen in general that $|D|$ is not pure dimensional, i.e., some maximal cells of $|D|$ can have dimension less than the dimension of $|D|$. For example if the tropical curve consists of four points $P,Q, R, S$ with an edge between $P$ and $Q$, two parallel edges from $P$ to $R$ and two parallel edge from $Q$ to $S$, and the divisor $D$ is equal to $P+Q$,  the complete linear system $|D|$  consists of three maximal cells: two segments and a triangle.

\subsection{The reduced divisor map}\label{sec:reduceddivisormap} 
Let $D$ be a divisor with $r(D) \geq 0$. For a point $P \in \Gamma$, we denote by $D_P$ the unique $P-$reduced divisor linearly equivalent to $D$. Note that $r(D) \geq 0$ is equivalent to $D_P \geq 0$ for every point $P$. Since all the coefficients of $D_P$ are already non-negative, this simply means that the coefficient  $D_P( P )$ of $D_P$ at $P$ is non-negative. 

\medskip 

Before stating the main theorem of this section, let us recall that a map $\phi$ from a tropical curve $\Gamma$ to a space $X$ with an abstract rational polyhedral complex structure $\Sigma$ is called integral affine if it is piece-wise linear with integral slopes in $X$. In other words, $\phi$ is integral affine if for any point $v$ of $\Gamma$ and any sufficiently small segment $I$ on an edge adjacent to $v$ whose image under $\phi$ entirely lies inside some $\tau$ in $\Sigma$, the image of $I$ in $|\tau|$ is a segment with slope in $N_{\tau,\mathbb Q}$.   

\medskip

The main result of this section is the following theorem, whose proof gives an explicit description of the behavior of the reduced divisor $D_P\sim D$ under infinitesimal changes of the base point $P$.

\begin{theorem}\label{thm:continuity}
Let $D$ be a divisor of non-negative rank $r(D) \geq 0$. The map $\Red: \Gamma \rightarrow |D| \hookrightarrow \Gamma^{(\deg(D))}$ defined by sending a point $P\in \Gamma$ to the point defined in $\Gamma^{\deg(D)}$   by the (unique) $P$-reduced  divisor $D_P$ linearly equivalent to $D$ is integral affine, and thus, continuous. 
\end{theorem}

\begin{proof}Let $P$ be a point of $\Gamma$ and $E=D_P$ be the unique $P$-reduced divisor linearly equivalent to $D$. We give an explicit behavior of $\Red$ around $P$ which shows in particular that $\Red$ is integral affine.

\noindent Let $\vec u$ be  a unit vector tangent to $\Gamma$ at $P$. (It is clear that there are $\deg(P)$ unit vectors tangent to $\Gamma$ at $P$.) We consider the interval $[P,P+\delta_0\vec u]$ for $\delta_0$ small enough, and show how the family of reduced divisors $D_v$ behave for $v$ a point of $\Gamma$ in this interval. 
We first chose $\delta_0$ small enough such that there is no point of the support of $E\, (=D_P)$ in the half-open interval $(P,\delta_0 \vec u\,]$. Since we will be dealing only with effective divisors, the first property of reduced divisors, namely the non-negativity of coefficients outside the base point is obviously verified for the divisor $E$ and for every point in this interval. Thus, we will consider below the second property concerning closed connected subsets, and show how the divisor $E$ has to be modified  for any point of the interval $[P,P+\delta_0\vec u]$ in order to respect this property. 

\medskip

The following cases can happen:

\vspace{.3cm}

\begin{itemize}
\item[(i)] $E$ itself is $Q$-reduced for all the points $Q \in [P,\delta_0 \vec u]$. (This case happens for example when the coefficient of $E$ at $P$ is zero, i.e., $E(P)=0$.)
\end{itemize}

\vspace{.3cm}

 In this case, the restriction of $\Red$ to $[P,P+\delta_0\vec u]$ is the constant map, i.e., $\Red$ contracts the interval $[P,P+\delta_0 \vec u]$ to a point $E$.  It is obvious that $\Red$ is integral affine restricted to this interval.

\vspace{.3cm}

\begin{itemize}
\item[(ii)] There are points $Q\in [P,\delta_0\vec u]$ for which $E$ is not $P$-reduced. (From the description we give below, it will be clear that in this case, $E$ cannot be reduced for any point of the interval $(P,\delta_0\vec u\,]$.)
 \end{itemize}

\vspace{.3cm}
This means that there exists a closed connected set $X$ which does not contain  $Q$ and such that for any point $v \in \partial X$, we have $E(v)\geq \deg_X^{out}(v)$. In this case $P$ should belong to $X$, otherwise the cut $X$ would violate the assumption that $E$ is $P$-reduced. By the choice of $\delta_0$, for all the points $v \in (P,P+\delta_0 \vec u]$, we have $E(v) = 0$. Since $X$ separates $P$ and $Q$ in the interval $[P,Q]$, we have $[P,P+\delta_0\vec u\,] \nsubseteq X$. We conclude that $X \cap [P,P+\delta_0 \vec u] = P$. Thus, we have a closed connected set $X$ with the following properties:
\begin{itemize}
\item[(1)] The intersection of $X$ with the interval $[P,P+\delta_0\vec u]$ contains only $P$.
\item[(2)] For all $v \in \partial X$, $E(v)\geq \deg_X^{out}(v)$.
\end{itemize}
In particular, this shows that  $E$ is not reduced for any point of the interval $(P,P+\delta_0\vec u\,]$

\medskip 

Let $\mathcal F$ be the family of all the closed connected subsets $X$ of $\Gamma$ with the above two properties. For any set $X \in \mathcal F$, and for any point $v \in \partial X$, we have $E(v) \geq \deg^{out}_X(v) \geq 1$, which implies that $\mathcal F$ is a finite set. 
Define $Y$ as the union of all the elements of $\mathcal F$, $Y := \bigcup_{X \in \mathcal F} X$. We will next show that $Y$ itself is a closed connected subset of $\Gamma$ which verifies Conditions (1) and (2) above. In other words, we will see that $Y$ is the maximal element of $\mathcal F$.

\medskip 

\noindent Property (1) easily follows by observing that $Y \cap [P,P+\delta_0\vec u] = \bigcup_{X \in \mathcal F}\, \bigl(\,X \cap [P,P+\delta_0\vec u\,]\,\bigr) = \{P\}.$  To show Property (2), let $v$ be a point of $\partial Y$. By the finiteness of $\mathcal F$,  there exists at least one $X \in \mathcal F$ such that $\deg^{out}_Y(v)\leq \deg^{out}_X(v)$. This clearly shows that $E(v) \geq \deg^{out}_Y(v)$.

\medskip 

Let $\delta_0>0$  be a small enough positive real such that for any point $v \in \partial Y$ and any unit vector $\vec w$ tangent to $\Gamma$ at $v$ and pointing outside $Y$, the whole interval $(v,v+\delta_0\vec w\,]$ lives outside $Y$ and does not contain any point of the support of $E$. 

Let  $\vec w_1,\dots,\vec w_s$ different from $\vec u$ be all the different unit vectors tangent to $\Gamma$ at points $v_1,\dots,v_s \in \partial Y$, respectively, and pointing out of $Y$. Remark that we might have $v_i=v_j$ for two different indices $i$ and $j$. The intervals $(v_i,v_i+\delta_0\vec w_i\,]$ together with $[P,P+\delta_0\vec u\,]$ form  (a segment on) all the out-going edges from $Y$.  

\medskip

\noindent To simplify the presentation below, define the excess of $E$ at $P$ with respect to the cut $Y$ to be the integer quantity $\mathrm{ex}_{E, Y}(  P ) :=  E{(}P{)}- \deg^{out}_Y( P )+1$. Note that  since $P$ is a boundary point of $Y$, and $Y$ verifies Property (2), we have $\mathrm{ex}_{E, Y}(  P ) \geq 1$.

 For any non-negative $\delta \leq \frac {\delta_0}{\mathrm{ex}_{E, Y}(  P )}$, the rational function $f_\delta:\Gamma \rightarrow \mathbb R$ is defined as follows. 

\begin{itemize}
\item The restriction of $f_\delta$ on $Y$ is zero;
\item The restriction of $f_\delta$ on any interval $[v_i,v_i+ \mathrm{ex}_{E, Y}(  P ) . \delta\,\vec w_i]$ is linear of slope $-1$ for any $1\leq i\leq s$;
\item The restriction of $f_\delta$ on the interval $[P,P+ \delta \vec u \,]$ is linear of slope $-\mathrm{ex}_{E, Y}(  P )$; and 
\item $f_\delta$ takes value $-\delta . \mathrm{ex}_{E, Y}(  P )$ at any other point of $\Gamma.$
\end{itemize}

\medskip 

 Obviously, by the properties of $Y$, the divisor $E+\div(f_\delta)$ is effective. For any $\delta \leq \frac{\delta_0}{\mathrm{ex}_{E, Y}(  P )}$ define the point $Q_\delta$ of $\Gamma$ by $Q_\delta=P+\delta\vec u$.  We claim that $E+\div(f_\delta)$ is $Q_\delta$-reduced.
This will in turn imply the theorem. Indeed, assuming the claim, it is fairly easy to check that $E+\div(f_\delta)$ is integral affine on the interval $[0,\frac{\delta_0}{\mathrm{ex}_{E, Y}(  P )}]$. More explicitly, the image  by the reduced divisor map $\Red$ of any point $Q_\delta $ in the interval $[\,P,P+\frac{\delta_0}{\mathrm{ex}_{E, Y}(  P ) }\vec u\,]$ will be equal to the point of $\Gamma^{(\deg(D))}$ whose coordinates are described as follows. 
Among the $\deg(D)$ points of $\Gamma$ in $\Red(Q_\delta)$, 
\begin{itemize}
\item $ \mathrm{ex}_{E, Y}(  P )$ of them are equal to $ Q_\delta = P+\delta \vec u$; 
\item $s $ of them are the points $v_1+\delta .\mathrm{ex}_{E, Y}(  P ) \vec w_i ,\dots,v_s + \delta.\mathrm{ex}_{E, Y}(  P )  \vec w_s$; 
\item and the remaining points are given by the image $ \Red\bigl(\sum_{v \neq P} (E(v) - \deg^{out}_Y(v))(v)\bigr)$ in $ \Gamma^{(\deg(D) -s -\mathrm{ex}_{E, Y}(  P ) )}$
\end{itemize}
(Here $\deg^{out}_Y(v)=0$ if $v$ is not a boundary point of $Y$.) Since $ \Red\bigl(\sum_{v \neq P} (E(v) - \deg^{out}_Y(v))(v)\bigr)$ does not depend on $\delta$ and so is a constant function on the interval  $[0,\frac{\delta_0}{\mathrm{ex}_{E, Y}(  P )}]$, and since the other points are integral affine functions of $\delta$, the whole map $\Red$ becomes integral affine on the interval $[P,P+ \frac{\delta_0}{\mathrm{ex}_{E, Y}(  P )}\vec u]$. This being true for any $P$ and any tangent vector $\vec u$ to $\Gamma$ at $P$, the statement of the theorem follows.  

\medskip 

 Thus, all we need to prove is to show that $E^\delta:=E+\div(f_\delta)$ is $Q_\delta$-reduced. For the sake of a contradiction, suppose this is not the case. Since $E^\delta$ is effective, this means there exists a closed connected subset $X$ of $\Gamma$ which does not contain $Q_\delta$ and such that for all $v$ on the boundary of $X$, we have $E^\delta(v) \geq \deg^{out}_X(v)$. First we claim that $X$ does not contain $P$, and more generally any point of the interval $[P,Q_\delta]$. Indeed, otherwise, since $X$ does not contain $Q_\delta,$ it should contain a point $v$ of $[P,Q_\delta)$ on its boundary. But clearly for this point we have $E^\delta(v) =0 < \deg^{out}_X(v)$ which is a contradiction by the choice of $X$. We now show that $X$ contains one of the points $v_i + \delta.\mathrm{ex}_{E, Y}(  P )\, \vec w_i$ on its boundary.  Indeed, since $E$ is $P$-reduced and $X$ is a closed and connected subset which does not contain $P$, by the definition of reduced divisors, there exists a point $v \in \partial X$ such that $E(v) < \deg^{out}_X(v)$. If $v \neq v_i + \delta. \mathrm{ex}_{E, Y}(  P )\,\vec w_i$ for all $1\leq i\leq s$, then by the construction of $E^\delta$, we have $E^\delta(v) \leq E(v) < \deg^{out}_X(v)$, which again leads to a contradiction by our choice of $X$. This shows our second claim. 

\medskip 

To summarize,  we have proved so far  that $X$ is a closed connected subset of $\Gamma$ which does not contain any point of $[P,Q_\delta]$, while $X$ contains at least one of the points $v_i+\delta.\mathrm{ex}_{E, Y}(  P ) \vec w_i$ on its boundary $\partial X$. Without loss of generality, let $v_1+\delta.\mathrm{ex}_{E, Y}(  P ) \vec w_1,\dots, v_k+\delta.\mathrm{ex}_{E, Y}(  P ) \vec w_k$, $k \geq 1$, be all the points among $v_i +\delta \vec w_i$ for  $i\leq s$ which lie on the boundary of $X$. 

\medskip 

For any point $v_i+\delta.\mathrm{ex}_{E, Y}(  P ) \vec w_i$ in $ \partial X$, we claim that either $[v_i,v_i+\delta.\mathrm{ex}_{E, Y}(  P ) \vec w_i]\subset X $ or $(v_i,v_i+\delta.\mathrm{ex}_{E, Y}(  P ) \vec w_i] \cap X = \{v_i+\delta.\mathrm{ex}_{E, Y}(  P ) \vec w_i\}$. To see this, note that  if none of these situations happen, then $X$ has to contain another point $x$ of $(v_i,v_i+\delta.\mathrm{ex}_{E, Y}(  P ) \vec w_i]$ on its boundary. For this point, obviously, we have $E^\delta(x) = 0 <\deg^{out}_X(x)$, which is certainly a contradiction by the choice of $X$. Thus, without loss of generality, assume $v_1+\delta.\mathrm{ex}_{E, Y}(  P ) \vec w_1,\dots, v_h+\delta.\mathrm{ex}_{E, Y}(  P ) \vec w_h$, $h\leq k$ are all the points on the boundary of $X$ with the property that 
$$(v_i, v_i+\delta.\mathrm{ex}_{E, Y}(  P ) \vec w_i] \cap X =\{v_i + \delta.\mathrm{ex}_{E, Y}(  P ) \vec w_i\}.$$ 

\medskip 

We divide the rest of the proof into two cases depending on whether $h=0$ or $h \geq 1$.

\medskip

\begin{itemize}
\item[$\bullet$] $h=0$; in other words, for any $i \leq k$, $[v_i,v_i+\delta .\mathrm{ex}_{E, Y}(  P )\vec w_i]\subset X$. 

Let $X' := X \setminus \bigcup_{i\leq k} (v_i,v_i+\delta.\mathrm{ex}_{E, Y}(  P ) \vec w_i\,]$. Clearly $X'$ is a closed connected subset of $\Gamma$ which does not contain $P$. 
Since $E$ is $P$-reduced, there exists a point $v \in \partial X'$ such that $E(v) < \deg^{out}_{X'}(v)$. 
We claim that there exists a $j\leq k$ such that $v=v_j$. Indeed,
if $v \neq v_j$ for all $j\leq k$, then $E^\delta(v) \leq E(v) < \deg^{out}_{X'}(v) = \deg^{out}_X(v)$, which clearly contradicts  the assumption we made on $X$. 
 Consider now the set $J$ of all indices $j \leq k$ such that $v_j =v$. Obviously, we have $\deg^{out}_{X'}(v) = \deg^{out}_X(v) + |J|$. By the definition of $E^\delta$, we also have $E^\delta(v) \leq E(v) -|J|$. It follows that $E^\delta(v) \leq E(v) - |J| < \deg^{out}_{X'}(v)-|J| = \deg^{out}_X(v)$.
Since $E^\delta$ is effective, we infer that $v \in \partial X$ and $E^\delta(v) < \deg^{out}_X(v)$. This is again in contradiction with our assumption on $X$. Thus, $E^\delta$ is $Q_\delta$-reduced in this case.

\item[$\bullet$]  $h \geq 1$.

In this case, define 
$$Y' = \Bigl ( Y \cup X \Bigr) \setminus \bigcup_{h+1 \leq i\leq k} (v_i,v_i + \delta.\mathrm{ex}_{E, Y}(  P ) \vec w_i\,].$$
Since $h \geq 1$, we obviously have $Y \subsetneq Y'$. Obviously, $Y'$ is closed and connected, and contains $P$. We will show that for all the points $v \in \partial Y'$, we have $E(v) \geq \deg^{out}_{Y'}(v)$, in other words, $Y' \in \mathcal F$ (for the family $\mathcal F$ defined at the beginning of the proof). This will be in contradictition with the maximality of $Y$ in $\mathcal F$, and thus, finishes the proof of the theorem.  Let $v \in \partial Y'$. Two cases can occur: 

\noindent Either, $v \in \{v_1,\dots, v_k, P\}$, in which case, by the properties of $Y$ and the fact that $Y \subset Y'$, we have 
$$E(v) \geq \deg^{out}_Y(v) \geq \deg^{out}_{Y'}(v).$$

\noindent Or, $v \in \partial X \setminus \{v_i + \delta.\mathrm{ex}_{E, Y}(  P )\vec w_i\}_{i=1}^k$ and $v$ is not equal to any of the points $P$, $Q_\delta$, and $v_i,$ for $i\leq k$. In this case, by the choice of $X$,  we have $E^\delta(v) \geq \deg^{out}_X(v)$.

\noindent By the definition of $E^\delta$, we have $E(v) = E^\delta(v)$. We again infer that 
$$E(v) \geq \deg^{out}_X(v) \geq \deg^{out}_{Y'}(v).$$ 
\end{itemize}
\end{proof}

\section{Ample Divisors and Canonical Embeddings}\label{sec:canonicalembeddings}
In the previous section, we proved that the map $\Red_D:\Gamma \rightarrow |D|$ is continuous. By analogy with the classical case, in this section we consider the cases where the map defines an embedding of the tropical curve. Since $|D|$ does not admit a tropical structure in general, here by an embedding we simply mean an injective (proper) map. The following definition could be thought of the tropical analogue  of ample and very ample divisors on algebraic curves, via embeddings into projective spaces (see~\cite{HMY09}).

\noindent Recall that $R(D)$ denotes the set of all rational functions $f$ on $\Gamma$ such that $D + \div(f) \geq 0$.

\begin{definition}[Ample and very ample divisors]\rm
A divisor $D$ is called very ample if the rational functions in $R(D)$ separates points of $\Gamma$, i.e., for all pair of points $P\neq Q$ in $\Gamma$, there exist two rational functions $f,g \in R(D)$ such that $f(P)-g(P) \neq f(Q)-g(Q)$. A divisor $D$ is called ample if an integer multiple of $D$ is very ample, i.e., there exists $m \in \mathbb N$ such that $m D $ is very ample.
\end{definition}

\begin{remark}\label{rem:ample}\rm
If $D$ itself is an effective divisor, i.e., if $0 \in R(D)$, the assertion $D$ is very ample is equivalent to the existence of a function $f\in R(D)$ such that $f(P) \neq f(Q)$.
\end{remark}

We have the following tropical analogue of the classical theorem on characterization of  very ample divisors in terms of embeddings into linear systems (for curves).
\begin{theorem}\label{thm:embedding}
A divisor $D$ is very ample if and only if the map $\Red$ defines an embedding of $\Gamma$ in $|D|$.
\end{theorem}
The proof is based on the following lemma.
\begin{lemma}\label{lem:maxreduced}
Let $P\in \Gamma$ and $D$ be a $P$-reduced divisor. A function $f \in R(D) $ takes its maximum value at $P$.
\end{lemma}
\begin{proof}
 For the sake of a contradiction suppose this is not the case and let $X$ be the set of all points $Q \in \Gamma$ where $f$ takes its maximum (thus, $P \notin X$). Note that $X$ is a closed subset of $\Gamma$. 
Since $f \in R(D)$, we have $D + \div(f) \geq 0$. Since $D$ is $P$-reduced and $X$ is closed and does not contain $P$, there is point $v \in \partial X$ such that $\deg^{out}_X (v) > D(v)$. Note that $f$ is strictly decreasing along any outgoing branch of $X$ at $v$, so the coefficient of $v$ in  $D + \div(f)$ is at most $D_{v} - \deg^{out}_X (v) < 0$. This contradicts the assumption that $D+\div(f)$ is effective and the lemma follows.
\end{proof}
\begin{proof}[Proof of theorem~\ref{thm:embedding}]
We show that $\Red_D$ is injective if and only if $D$ is very ample. 
\begin{itemize}
\item Let $D$ be very ample. For the sake of a contradiction, suppose that $\Red_D$ is not injective. This means there are two points $P$ and $Q$ such that the $P$- and $Q$-reduced divisors linearly equivalent to $D$ are the same, i.e., $D_P = D_Q$. For any rational function $f \in R(D_P)=R(D_Q)$, by applying Lemma~\ref{lem:maxreduced}, we know that $f$ takes its maximum value both at $P$ and $Q$. This means that $f(P) = f(Q)$ for all $f \in R(D_P)$ and this contradicts the assumption that $D$ is a very ample divisor (see Remark~\ref{rem:ample}).
\item If $\Red_D$ is injective, then $D_P \neq D_Q$. For two points $P\neq Q$ in $\Gamma$, there exists a non-constant rational function $g \in R(D_P)$ such that $D_P + \div(g)=D_Q$. By Lemma~\ref{lem:maxreduced} applied twice, $g$ takes its maximum value at $P$ and $-g$ takes its maximum value at $Q$. Since $g$ is not constant, this means that $g(P)\neq g(Q)$ and so $D$ is very ample. 
\end{itemize}
\end{proof}

\begin{corollary}\label{thm:2g+1} If $\deg(D) \geq 2g+1$, then $D$ is very ample. In particular, a divisor $D$ is ample if and only if $\deg(D) >0$.
\end{corollary}

\begin{proof}
This is a consequence of Theorem~\ref{thm:embedding}. Let $E=D_P \sim D$ be $P$-reduced. By Riemann-Roch theorem, $E$ is effective. We claim that the coefficient of $P$ in $E$ is at least $\deg(D)-g$. Consider the set $S$ obtained as the union of all the branching points of $\Gamma$ and all the points in the support of $E$ 
$$S := \mathrm{support}\,(E) \cup \{\textrm{ branching points of }\Gamma\,\}.$$
Suppose $|S|=N+1$ and define greedily an ordering $v_0, \dots, v_N$ of $S$ as follows. Start by setting $v_0= P$. Suppose by induction that the partial ordering $v_0,\dots,v_i$ is already defined. Consider $X_i$ the cut defined by all the points of $S \setminus\{v_0,\dots,v_i\}$. This is basically obtained by taking all the segments joining two points of $S \setminus\{v_0,\dots,v_i\}$ in $\Gamma$. Since $P \notin X_i$ and $E$ is $P$-reduced, there exists a point $v$ such that $E(v) < \deg^{out}_{X_i}(v)$. Define $v_{i+1}=v$. It is straightforward to check that that $\sum_{i\geq 1} \deg^{out}_{X_i}(v_{i+1}) -1 = g$. Thus, by the construction of the ordering, we have

\[\sum_{i\geq 1} E(v_i) \leq \sum_{i\geq 1}\deg^{out}_{X_i}(v_{i+1}) -1 = g.\]
This clearly implies that $E( P )=E(v_0) = \deg(D) - \sum_{i\geq 1} E(v_i) \geq \deg(D) -g$, which establishes the claim.

To conclude the proof, observe that if $\deg(D) \geq 2g+1$, it cannot happen that $D_P = D_Q$. Otherwise, we would have 
$$\deg(D)\geq D_P( P )+D_P(Q) = D_P( P ) + D_Q( Q )\geq \deg(D)-g +\deg(D)-g \geq \deg(D)+1.$$ 
By Theorem~\ref{thm:embedding}, $D$ is very ample provided that $\deg(D)\geq 2g+1$.  
\end{proof}

\medskip 
We end this section by giving an improvement on the statement of Corollary~\ref{thm:2g+1}. 
First, we note that the bound $\deg(D) \geq 2g+1$ in the above corollary is tight. To give an example, consider a tropical curve $\Gamma_g$ obtained by associating positive lengths to the edges of the banana graph of genus $g$ with two branching points $P$ and $Q$ of degree $g+1$ connected by $g+1$ parallel edges, see Figure~\ref{fig:case1}. The divisor $D=g\,( P )+g\,( Q )$ on $\Gamma_g$ has degree $2g$ and is both $P$- and $Q$-reduced, so $D$ is not very ample. However, this is the only situation where this happens. More precisely, if  $\Gamma$ does not have any model whose underlying graph is isomorphic to the banana graph of genus $g$, i.e., the graph underlying $\Gamma_g$, then every divisor $D$ of degree at least $2g$ is very ample. We show this by absurd, so for the sake of a contradiction, consider a divisor $E$ of degree $2g$ which is both $P$- and $Q$-reduced for two different points $P \neq Q$ of $\Gamma$. As the proof of Corollary~\ref{thm:2g+1} shows, the coefficient of $P$ and $Q$ in $E$ are at least $g$, and since $\deg(E)=2g$ and $E$ is effective, then we must have $E = g\,( P )+g\,( Q )$. Since $E$ is $P$-reduced, by taking the cut consisting of  the only point $Q$, we see that $\deg_\Gamma(Q) \geq g+1$. A similar argument shows that $Q$ is not a cut vertex, i.e.,  $\Gamma \setminus \{Q\}$ is connected. Similarly $\deg_\Gamma(P) \geq g+1$ and $P$ is not a cut-vertex. The genus of $\Gamma$ being $g$, we infer that $\deg_\Gamma(P)=\deg_\Gamma(Q) = g+1$ and all the other points are of degree two. Since none of $P$ and $Q$ is a cut-vertex, this shows that the underlying graph of $\Gamma$ consists of $P$ and $Q$ and $g+1$ parallel edges between them. 

\subsection{Canonical embedding} In this section, we consider the map $\Red_K$ defined by the canonical divisor $K=K_\Gamma$ of a tropical curve $\Gamma$. We provide a complete characterization of all the situations where the canonical divisor is not very ample. 

\medskip

\noindent We start this section by proving the following simple lemma.

\begin{lemma}\label{lem:canonic}
Let $\Gamma$ be a tropical curve of genus $g \geq 2$. If $K$ is not very ample, then there are two points $P$ and $Q$ of $\Gamma$ such that $K_P=K_Q = (g-1).(P)+(g-1).(Q)$. Here $K_P $ and $ K_Q$ are the $P$- and $Q$-reduced divisors linearly equivalent to $K$.
\end{lemma}

\begin{proof}
We again use Theorem~\ref{thm:embedding}. If $K$ is not very ample, then the map $\Red_K$ is not injective, i.e., there exist two points $P$ and $Q$ such that $K_P =K_Q$. 
We claim that the coefficient of $P$ in $K_P$ is at least $g-1$.  Indeed, by Riemann-Roch theorem, $K - (g-1)(P)$ has non-negative rank. This means that the $P$-reduced divisor associated to $K - (g-1)(P)$ is effective. This $P$-reduced divisor being simply $K_P - (g-1)( P )$, we infer that the coefficient of $P$ in $K_P$ is at least $g-1$. The same argument gives $K_P(Q)=K_Q(Q)\geq g-1$. Since $\deg(K) = 2g-2$, and all the other coefficients of $K_P$ are non-negative, we must have $K_P = K_Q= (g-1).(P) + (g-1).(Q)$.
\end{proof}

\begin{figure}[!htb] 
  \begin{minipage}[t]{.45\textwidth}
    \begin{center}  
  \subfigure[\label{fig:case1}]{\includegraphics[width=0.4\linewidth]{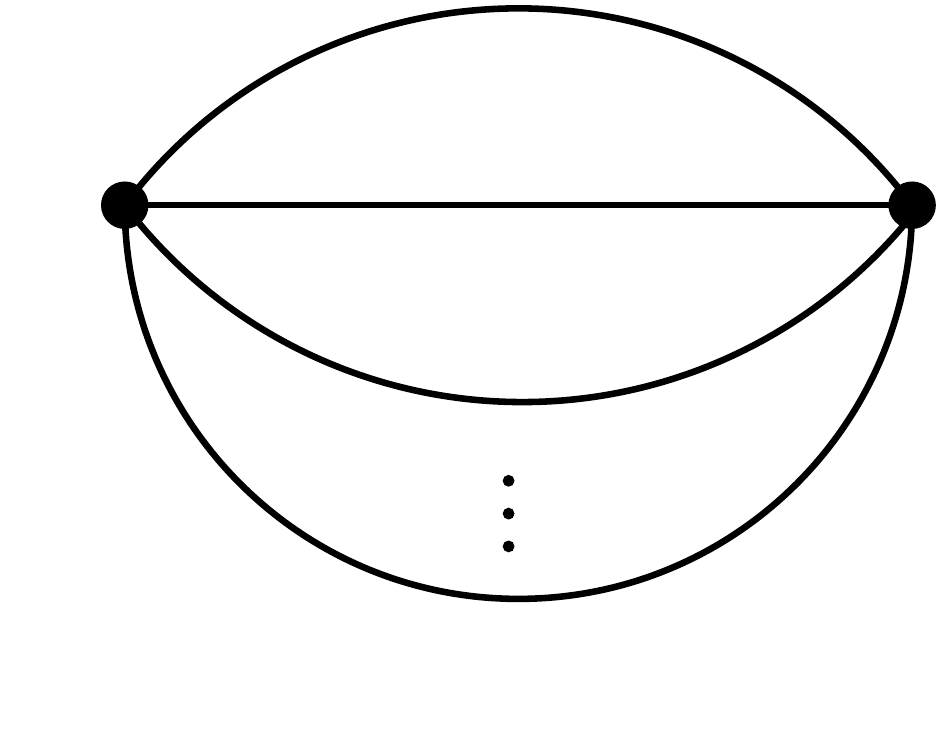}}  
  \hspace{1cm}
    \subfigure[\label{fig:case2}]{\includegraphics[width=0.4\linewidth]{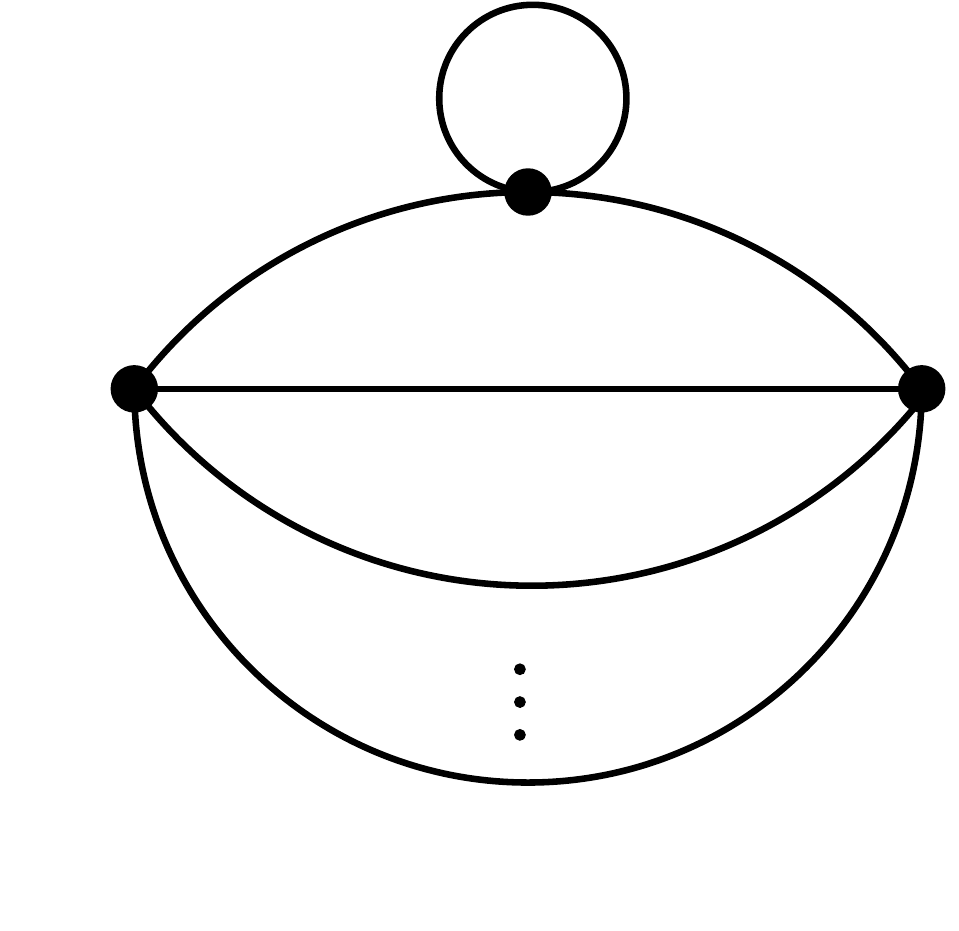}}
        \end{center}
  \end{minipage}
  \hspace{.7cm}
  \begin{minipage}[t]{.45\textwidth}
    \begin{center}  
  \subfigure[\label{fig:case3}]{\includegraphics[width=0.4\linewidth]{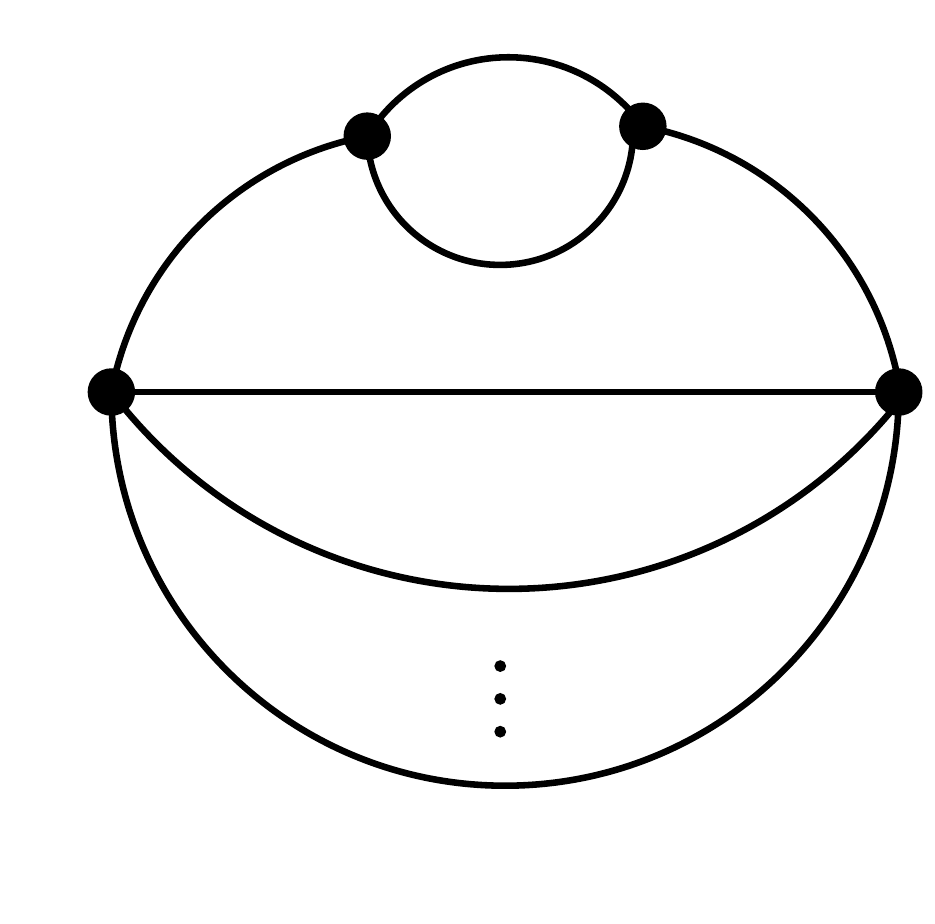}}
     \hspace{.8cm}
     \subfigure[\label{fig:case4}]{\includegraphics[width=0.4\linewidth]{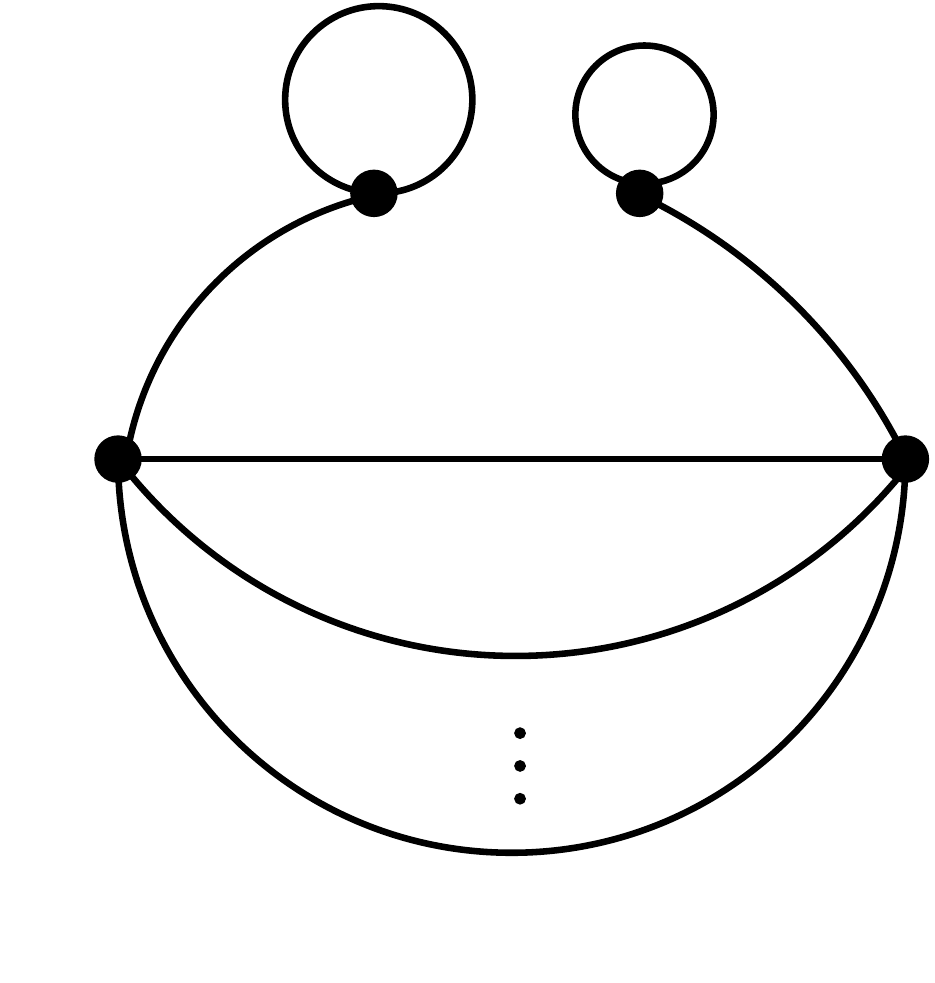}}

   \setlength{\unitlength}{4144sp}%
\begingroup\makeatletter\ifx\SetFigFont\undefined%
\gdef\SetFigFont#1#2#3#4#5{%
  \reset@font\fontsize{#1}{#2pt}%
  \fontfamily{#3}\fontseries{#4}\fontshape{#5}%
  \selectfont}%
\fi\endgroup%
\begin{picture}(4260,4486)(4276,-6182)
\put(700,-700){\makebox(0,0)[lb]{\smash{{\SetFigFont{10}{20.0}{\familydefault}{\mddefault}{\updefault}{\color[rgb]{0,0,0}$P$}%
}}}}
\put(2000,-710){\makebox(0,0)[lb]{\smash{{\SetFigFont{10}{20.0}{\familydefault}{\mddefault}{\updefault}{\color[rgb]{0,0,0}$Q$}%
}}}}

\put(2550,-700){\makebox(0,0)[lb]{\smash{{\SetFigFont{10}{20.0}{\familydefault}{\mddefault}{\updefault}{\color[rgb]{0,0,0}$P$}%
}}}}
\put(3820,-710){\makebox(0,0)[lb]{\smash{{\SetFigFont{10}{20.0}{\familydefault}{\mddefault}{\updefault}{\color[rgb]{0,0,0}$Q$}%
}}}}

\put(4320,-700){\makebox(0,0)[lb]{\smash{{\SetFigFont{10}{20.0}{\familydefault}{\mddefault}{\updefault}{\color[rgb]{0,0,0}$P$}%
}}}}
\put(5615,-710){\makebox(0,0)[lb]{\smash{{\SetFigFont{10}{20.0}{\familydefault}{\mddefault}{\updefault}{\color[rgb]{0,0,0}$Q$}%
}}}}

\put(6070,-700){\makebox(0,0)[lb]{\smash{{\SetFigFont{10}{20.0}{\familydefault}{\mddefault}{\updefault}{\color[rgb]{0,0,0}$P$}%
}}}}
\put(7360,-710){\makebox(0,0)[lb]{\smash{{\SetFigFont{10}{20.0}{\familydefault}{\mddefault}{\updefault}{\color[rgb]{0,0,0}$Q$}%
}}}}

\put(6520,-480){\makebox(0,0)[lb]{\smash{{\SetFigFont{10}{20.0}{\familydefault}{\mddefault}{\updefault}{\color[rgb]{0,0,0}$R$}%
}}}}

\put(6900,-480){\makebox(0,0)[lb]{\smash{{\SetFigFont{10}{20.0}{\familydefault}{\mddefault}{\updefault}{\color[rgb]{0,0,0}$S$}%
}}}}

\put(4680,-250){\makebox(0,0)[lb]{\smash{{\SetFigFont{10}{20.0}{\familydefault}{\mddefault}{\updefault}{\color[rgb]{0,0,0}$R$}%
}}}}

\put(5270,-250){\makebox(0,0)[lb]{\smash{{\SetFigFont{10}{20.0}{\familydefault}{\mddefault}{\updefault}{\color[rgb]{0,0,0}$S$}%
}}}}

\put(3210,-530){\makebox(0,0)[lb]{\smash{{\SetFigFont{10}{20.0}{\familydefault}{\mddefault}{\updefault}{\color[rgb]{0,0,0}$R$}%
}}}}
\end{picture}%
    \end{center}
    \end{minipage}
\vspace{-10cm}
    \caption{Cases where $K$ is not very ample, c.f.,~Theorem~\ref{thm:nottropical}.\newline (a) Case C.$\mathrm{I}$: $\deg(P)=\deg(Q)=g+1$, the lengths are arbitrary.
      \newline
      (b) Case C.$\mathrm{II}$: $\deg(P)=\deg(Q)=g, \ \deg(R)=4$, and $\ell([P,R])=\ell([R,Q])$. The other lengths are arbitrary.\newline
      (c) Case C.$\mathrm{II'}$: $\deg(P)=\deg(Q)=g, \ \deg(R)= \deg(S)=3$, and $\ell([P,R])=\ell([S,Q])$. The other lengths are arbitrary.\newline
      (d) Case C.$\mathrm{III}$: $\deg(P)=\deg(Q)=g, \ \deg(R)= \deg(S)=3$, and the lengths are arbitrary.
      }

\end{figure}

The following theorem characterizes all the situations where the canonical divisor is not very ample. We  can assume that $\Gamma$ does not  have any leaf.
\begin{theorem}\label{thm:nottropical}
Let $\Gamma$ be a tropical curve of genus $g \geq 2$ without points of degree one. The only cases where the canonical divisor is not very ample are the following:
\begin{itemize}
\item[$\mathrm{C}.\mathrm{I}$] The underlying graph of $\Gamma$ consists of  two vertices $P$ and $Q$ of degree $g+1$ which are connected by $g+1$ parallel edges. In particular,  $K = (g-1)(P)+(g-1)(Q)$.
\item[$\mathrm C.\mathrm{II}$] The underlying graph of $\Gamma$ consists of two vertices $P$ and $Q$ of degree $g$ which are joined by $g$ parallel edges and a point of degree four in the middle of one of the edges with an extra edge from this point to itself, i.e., a point $R$ such that the segment $[P,R]$ and $[R,Q]$ have the same length and there is a loop connecting $R$ to $R$ (see \emph{Figure~\ref{fig:case2}}).
\item[$\mathrm C.\mathrm{II'}$] The underlying graph of $\Gamma$ consists of four branching points $P$, $Q$, $R$, and $S$ as in \emph{Figure~\ref{fig:case3}}, with $\deg(P)=\deg(Q)=g$ and $\deg(R)=\deg(S) =3$. In addition $[P,R]$ and $[Q,S]$ have the same length.
\item[$\mathrm C.\mathrm{III}$]  The underlying graph of $\Gamma$ have four branching points $P$, $Q$, $R$, and $S$ as in \emph{Figure~\ref{fig:case4}}, with $\deg(P)=\deg(Q)=g$ and $\deg(R)=\deg(S) =3$. The lengths are arbitrary.
\end{itemize} 
\end{theorem}

Recall that a tropical curve is called {\it hyperelliptic} if there exists a divisor $D$ of degree two  such that $r(D) =1$. 
For a tropical curve whose canonical divisor is not very ample, let $D = (P)+(Q)$, such that $K\sim (g-1).D$ by Lemma~\ref{lem:canonic}. By the above characterization theorem, it is easy to see that $r(D)=1$. Since $\deg(D)=2$, the tropical curves $\Gamma$ appearing in the characterization theorem above are all hyperelliptic. In addition, one obtains that canonical divisors of generic tropical curves of genus $g \geq 2$ are very ample.

\begin{figure}
  \begin{center}
  \subfigure[\label{fig:impossible1}]{\includegraphics[width=0.2\linewidth]{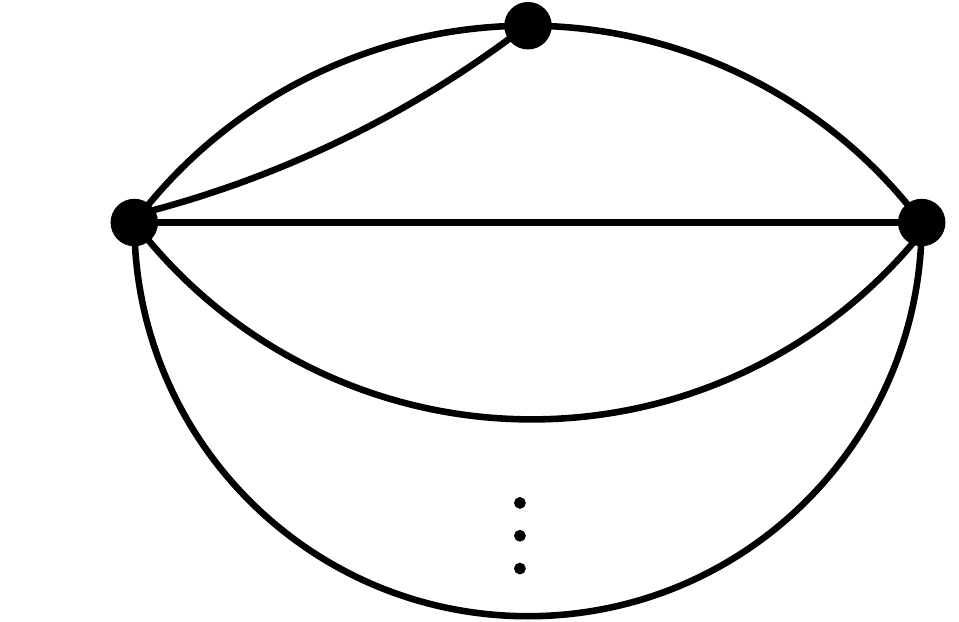}} \hspace{4cm}
  \subfigure[\label{fig:impossible2}]{\includegraphics[width=0.2\linewidth]{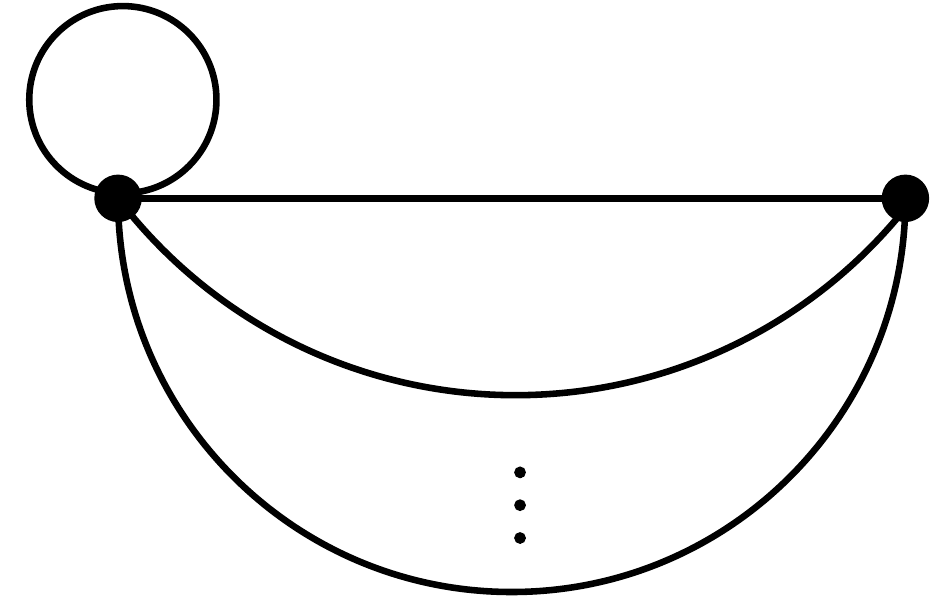}}
\end{center}
 \setlength{\unitlength}{4144sp}%
\begingroup\makeatletter\ifx\SetFigFont\undefined%
\gdef\SetFigFont#1#2#3#4#5{%
  \reset@font\fontsize{#1}{#2pt}%
  \fontfamily{#3}\fontseries{#4}\fontshape{#5}%
  \selectfont}%
\fi\endgroup%
\begin{picture}(4260,4486)(3276,-6000)
\put(4400,-600){\makebox(0,0)[lb]{\smash{{\SetFigFont{10}{20.0}{\familydefault}{\mddefault}{\updefault}{\color[rgb]{0,0,0}$P$}%
}}}}
\put(5850,-600){\makebox(0,0)[lb]{\smash{{\SetFigFont{10}{20.0}{\familydefault}{\mddefault}{\updefault}{\color[rgb]{0,0,0}$Q$}%
}}}}

\put(5150,-200){\makebox(0,0)[lb]{\smash{{\SetFigFont{10}{20.0}{\familydefault}{\mddefault}{\updefault}{\color[rgb]{0,0,0}$R$}%
}}}}

\put(7850,-450){\makebox(0,0)[lb]{\smash{{\SetFigFont{10}{20.0}{\familydefault}{\mddefault}{\updefault}{\color[rgb]{0,0,0}$P$}%
}}}}
\put(9050,-450){\makebox(0,0)[lb]{\smash{{\SetFigFont{10}{20.0}{\familydefault}{\mddefault}{\updefault}{\color[rgb]{0,0,0}$Q$}%
}}}}
\end{picture}
\vspace{-10cm}
\caption{(a) $\deg(P)=g+1$, $\deg(Q)=g$, $\deg(R)=3$.\newline
   (b) $\deg(P)=g+2$, $\deg(Q)=g$.}
\end{figure}

\begin{proof}[Proof of Theorem~\ref{thm:nottropical}]
Suppose $K$ is not very ample. By Lemma~\ref{lem:canonic}, there are two points $P$ and $Q$ such that $K_P = K_Q = (g-1)(P)+(g-1)(Q)$.  Since $K_P$ is $P$-reduced, we have $\deg(Q) \geq g$. Similarly we have $\deg(P) \geq g$. Since $\sum_{P \in \Gamma} \deg(P)-2 = 2g-2$ and $\Gamma$ does not have any leaf, the following cases happen:
\begin{itemize}
\item $\deg(P) = \deg(Q) =g+1$. In this case, $K_P=K_Q$ being $P$-(resp. $Q$-)reduced, $Q$ (resp. $P$) cannot be a cut-vertex in $\Gamma$. Thus, the underlying graph of $\Gamma$ consists of  two vertices $P$ and $Q$ of degree $g+1$ which are connected by $g+1$ parallel edges. In particular,  $K = (g-1)(P)+(g-1)(Q)$. And this is case $C.\mathrm I$ of the theorem.
\item  Either $\deg(P) = g+1$ and $\deg(Q) =g$, or $\deg(P) = g$ and $\deg(Q) =g+1$. Without loss of generality, suppose that $\deg(P) = g+1$ and $\deg(Q) =g$. In this case, $\Gamma$ has a branching point $R\neq P,Q$ of degree three. The underlying graph of $\Gamma$ is the graph of Figure~\ref{fig:impossible1}. It is easy to see that since $R \neq Q$, $K \nsim  (g-1)(P)+(g-1)(Q)$, so this situation cannot happen.

\item $\deg(P) = \deg(Q) =g$. In this case, $\Gamma$ has either one branching point $R \neq P, Q$ of degree four, or two different branching points $R\neq P,Q$ and $S\neq P,Q$ of degree three. The underlying graph of $\Gamma$ has then a very particular shape. A case analysis shows that the only possibilities for $(g-1)(P)+(g-1)(Q)$ to be equivalent to the canonical divisor are the  cases $C.\mathrm{II}$, $C.\mathrm{II'}$, and $C.\mathrm{III}$ described in the theorem, see Figures~\ref{fig:case2},\ref{fig:case3},\ref{fig:case4}.
\item $\deg(P)=g$ and $\deg(Q) =g+2$ (or  $\deg(P)=g+2$ and $\deg(Q) =g$). In this case, the underlying graph of $\Gamma$ consists of the two points $P$ and $Q$ joined by $g$ parallel edges and there is an extra loop based on $P$ (see Figure~\ref{fig:impossible2}). It is easy to show that this case cannot happen. 
\end{itemize}
 The theorem follows.
\end{proof}

\section{Applications: Tropical Weierstrass Points and Rank-Determining Sets}\label{sec:weier}
We already saw some applications of the reduced divisor map, specially in obtaining a characterization of all the tropical curves with a very ample canonical divisor. In this section, we provide two more applications of this method: we first provide an effective way of proving the existence of Weierstrass points on tropical curves of genus at least two, and then  proceed by giving an alternative shorter proof of a theorem of Luo on tropical rank-determining sets.

\subsection{Tropical Weierstrass points}
 Let $\Gamma$ be a tropical curve of genus $g$, and $K$ the canonical divisor of $\Gamma$. By analogy with the theory of algebraic curves, we say that $P \in \Gamma$ is a Weierstrass point if $r(g ( P )) \geq 1$. 
Equivalently, by  Riemann-Roch theorem, a point $P \in \Gamma$ is Weierstrass if and only if $r(K- g(P)) \geq 0$, i.e., if $K- g(P)$ is equivalent to an effective divisor. Riemann-Roch theorem also implies

\begin{proposition}\label{prop:basic}
For every point $P\in \Gamma$,  $r(K -(g-1)(P)) \geq 0$. Let $K_P = a_P (P) + \sum_{v : v \neq p}a_v(v)$ be the (unique) $P$-reduced divisor linearly equivalent to $K$. Then $a_P \geq g-1$. In addition, $P$ is Weierstrass if and only if $a_P \geq g$. 
\end{proposition}

 We show in this section the following tropical analogue of the classical theorem that every smooth projective curve of genus at least two over an algebraically closed field $\kappa$ has a Weierstrass point. 

\begin{theorem}\label{thm:main}
 Any tropical curve of genus at least two has a Weierstrass point. 
 \end{theorem}
Before going through  the proof, the following remarks are in order.
\begin{remark}\label{rem:weier}\rm 
\begin{itemize}
\item The theorem easily follows  by Baker's specialization lemma~\cite{Bak07}, by using the existence of a degeneration of  smooth curves to $\Gamma$, from the classical theorem on the existence of Weierstrass points, see~\cite{Bak07} for more details. The proof given here has the advantage of being more explicit, providing an efficient way to finding a Weierstrass point.

\item In the classical setting, it is possible to count the number of Weierstrass points in a very precise sense, see Section~\ref{app:weier}. In the tropical context, such a formula cannot exist, at least with the definition of Weierstrass points given above. Indeed, as was pointed out in~\cite{Bak07}, the example of the graph $\Gamma_g$ with two branching points of degree $g+1$ with $g+1$ parallel edges for $g\geq 3$ (Figure~\ref{fig:case1})  shows that $\Gamma$ can have an infinite number of such points. In this tropical curve,  it is quite easy to see that there exists an interval of positive length on each edge joining $P$ to $Q$ whose  points are all Weierstrass.
\end{itemize}
\end{remark}

\begin{proof}[Proof of Theorem~\ref{thm:main}]Let $\Gamma$  be a tropical curve of genus $g \geq 2$ and $K$ the canonical divisor of $\Gamma$. Let $\Red : \Gamma \rightarrow \Gamma^{(2g-2)}$ be the map defined in Section~\ref{sec:reducedmap} which sends a point $P \in \Gamma$ to the point of $\Gamma^{(2g-2)}$ which corresponds to the unique $P$-reduced divisor $K_P$ linearly equivalent to $K$.
Our proof will be essentially based on Theorem~\ref{thm:continuity} and the explicit description of the behavior of the map $\Red$ given in the proof of Theorem~\ref{thm:continuity}.

\noindent Suppose for the sake of a contradiction that there is no Weierstrass point in $\Gamma$. Then, by Proposition~\ref{prop:basic}, for every point $P \in \Gamma$, the coefficient of $P$ in $K_P$ is exactly $g-1$. Let $F:\Gamma \rightarrow \mathbb R_{\geq 0}$ be the continuous function defined as follows. For a point $P \in \Gamma$, let $P_1,\dots,P_{g-1}$ be all the points in the support of $K_P$ which are different from $P$ (potentially we could have $P_i = P_j$ for some $i \neq j$). The value of $F$ at $P$ is defined as $F(P):=\min \Bigl\{\, \dist_\Gamma(P,P_i)\,\Bigr\}_{i=1}^{g-1}$. Recall that $\dist_\Gamma(.\,,.)$ denotes the distance between points in $\Gamma$. By Theorem~\ref{thm:continuity}, the point $(P_1,\dots,P_{g-1}) \in \Gamma^{(g-1)}$ is a continuous function of $P$, thus, the map $F$ is continuous.  By the compactness of $\Gamma$ and continuity of $F$, there is a point $P \in \Gamma$ such that $F$ takes its minimum value at $P$, i.e., $F(P)$ is the minimum of $F(Q)$ over $Q\in \Gamma$. By our assumption, all the points $P_i$ are different from $P$, and so $F(P) > 0$. Let $P_1\neq P$ be the point in the support of $K_P$ such that $F(P) = \dist_\Gamma(P,P_1)$ and let $C$ be a shortest path from $P$ to $P_1$. We show that moving $P$ to another point $P'$ in the direction of the shortest path $C$ will decrease the value of $F$, which will obviously be a contradiction.  

\noindent Let $\vec u$ be the unit tangent vector to $\Gamma$ at $P$ such that $[P,P+\delta \vec u\,] \subset C$ for all sufficiently small $\delta$. We look at the behavior of $K_P$ in the interval $[P,P+\delta\vec u\,]$. Since for a point $Q$ in this interval, the coefficient of $K_Q$ at $Q$ is at least $g-1 \geq 1$,  $K_P$ is evidently not $Q$-reduced for $\delta$ small enough. Thus, let $Y$ be the maximal closed connected subset of $\Gamma$, defined in the proof of Theorem~\ref{thm:continuity} such that 
\begin{itemize}
\item[1.] The intersection of $Y$ with the interval $[P,P+\delta \vec u]$ contains only $P$.
\item[2.] For all $v \in \partial Y$, $K_{P}(v) \geq \deg^{out}_Y(v)$.
\end{itemize}  

Note that, by our assumption, we also have  $K_{P}( P ) = g-1$.  The proof of Theorem~\ref{thm:continuity} shows that for $\delta$ small enough, and for $Q \in [P,P+\delta\vec u\,]$, the coefficient of $K_{Q}$ at $Q$ is equal to $K_{P}( P ) - \deg^{out}_Y(P)+1$.  Since $K_{P}( P )=K_{Q}(Q) =g-1$, we must have $\deg^{out}_Y(P) = 1$. 

\vspace{.3cm}
Two cases can happen: 
\begin{itemize}
\item[$\bullet$] Either $P_1 \notin Y$, in which case $K_{P}(P_1)=K_{Q}(P_1)$. In other words, $P_1$ lies in the support of $K_Q$ for all $Q \in [P,P+\delta\vec u\,]$ for $\delta>0$ small enough. However, note that $F(Q) \leq \dist_\Gamma (Q,P_1) < \mathrm{length}( C )=F(P)$, which leads to a contradiction. 

\item[$\bullet$] Or, $P_1 \in Y$. We claim that in this case $P_1 \in \partial Y$ and $C \cap Y = \{P_1,P\}$. Suppose for the sake of a contradiction that this is not true, i.e., $P_1$ lies in the interior of $Y$ or $C$ intersects $Y$ in a point different from $P_1$ and $P$. Since $C$ is a path from $P$ to $P_1$ and since the interval $[P,P+\delta\vec u\,] \subset C$ and $(P,P+\delta\vec u\,] \cap Y =\emptyset$, in either of these cases, there must  exist a point $v \in \partial Y$, $v\neq P_1, P$, which is the first time $C$ returns to $Y$ while going from $P$ to $P_1$. By Property 2. of the set $Y$ above, we have $K_{P}(v) \geq \deg^{out}_Y(v) \geq 1$, i.e., $v$ lies in the support of $K_P$. We have $\dist_\Gamma(P,v) < \dist_\Gamma(P,P_1) = F(P)$ which contradicts the definition of $F$. This proves our claim.

\noindent Let $\vec w$ be the unit vector tangent to $\Gamma$ at $P_1$ such that the interval $[P_1,P_1+\delta \vec w\,]$ lies inside $C$ for $\delta$ small enough. By our explicit description of the reduced divisors $K_Q$ for $Q \in [P,P+\delta\vec u\,] $, given in the proof of Theorem~\ref{thm:continuity}, there is a point $Q_1 \in [P_1,P_1+\delta \vec w\,]$  such that the coefficient of $Q_1$ in $K_Q$ is one, provided that $Q$ is close enough to $P$. This shows that $F(Q) \leq \dist_\Gamma(Q,Q_1) < \mathrm{length}( C ) = F(P)$, which again leads to a contradiction by our choice of the point $P$.
\end{itemize}

The proof of Theorem~\ref{thm:main} is now complete. 
\end{proof}

 The exact same argument gives the following more general statement.

\begin{theorem}\label{thm:weiergen} Let $D$ be a divisor of positive rank. If for all $P\in \Gamma$, the coefficient of $D_P$ at $P$ is at least $a$, for some integer $a\geq 1$, then there exists a point $P \in \Gamma$ such that the coefficient of $P$ in $D_P$ is at least $a+1$.
\end{theorem}

\begin{definition} \rm Let $D$ be a divisor of rank $r$. A point $P\in \Gamma$ is called a Weierstrass point for $D$ (or a $D$-Weierstrass point) if the coefficient of $D_P$ at $P$ is at least $r+1$.  
\end{definition}

\begin{corollary}
Let  $\Gamma$ be a tropical curve and $D$ be a divisor of positive rank. Then $\Gamma$ has a $D$-Weierstrass point.  
\end{corollary}
\begin{proof} Let $r = r(D) \geq 1$. By definition of rank, for any point $Q\in \Gamma$, the coefficient of $Q$ in $D_Q$ is at least $r$. By Theorem~\ref{thm:weiergen}, there exists a point $P \in \Gamma$, with $D_{P}( P ) \geq r+1$, i.e., $P$ is a $D$-Weierstrass point. 
\end{proof}

It is worth noticing that Remark~\ref{rem:weier} remains still valid for $D$-Weierstrass points. It would be interesting to have a definition of tropical Wronskians (as in the classical setting, c.f., Section~\ref{app:weier}) which would allow to obtain more precise statements about the number of such points.

\subsection{Rank-determining sets of points} Recall that a subset $A \subset \Gamma$ is called rank-determining~\cite{Luo09} if the rank of a divisor in $\Gamma$ is equal to the maximum $r$ such that for any effective divisor $E$ of degree $r$ with support in $A$, the divisor $D-E$ has non-negative rank (equivalently, $D-E$ is equivalent to an effective divisor).

\medskip

The following proposition can be easily proved by induction on $r$. 

\begin{proposition}\label{prop:rank-determ}
A given subset $A \subset \Gamma$ is rank-determining  if and only if for a divisor $D$ of non-negative rank the following two assertions are equivalent. 
\begin{itemize}
\item Rank of $D$ is at least one, $r(D) \geq 1$.
\item For any point $P \in A$, $P$ is in the support of $D_P$. 
\end{itemize}
\end{proposition}

Let $A$ be a set of points of $\Gamma$. By cutting $\Gamma$ along $A$ we mean the (possibly disconnected) tropical curve obtained as the disjoint union of  the closure in $\Gamma$ of each of the connected components of $\Gamma \setminus A$.

\begin{theorem}[Luo~\cite{Luo09}]\label{thm:f-width}
Given a subset $A \subset \Gamma$, if cutting $\Gamma$ along the points of $A$ results in a disjoint union of  at least two tropical curves, and all the connected components of the cutting are of genus zero, then $A$ is rank-determining. 
\end{theorem}
Before giving the proof, the following remarks are in order. 

\begin{remark}
\begin{itemize}

\item[$(1)$] For a model $G=(V,E)$ of a tropical curve $\Gamma$ without loops, the vertices of $V$ form a rank-determining set. This was essentially proved in~\cite{HKN08}.
\item[$(2)$] Luo's theorem is a more general statement in terms of special open sets (see~\cite{Luo09} for the definition), and gives a necessary and sufficient condition for a subset $A$ of $\Gamma$ to be rank-determining. We restrict ourselves to the above statement, however, we note that the proof given below works without any change in that more general setting and derives the sufficiency of the condition in Luo's theorem.

\end{itemize}
\end{remark}
\begin{proof}[Proof of Theorem~\ref{thm:f-width}]
Using  Proposition~\ref{prop:rank-determ} (twice, once for $A$ and once for $\Gamma$), we need to prove that if $P$ belongs to the support of $D_P$ for every point $P$ in $A$, and cutting $\Gamma$ along $A$ results in a set of at least two tropical rational curves (i.e., of genus zero), then  for any point $Q \in \Gamma$, $Q$ is in the support of $D_Q$. Let $\Gamma_1,\dots,\Gamma_s$ be all the tropical curves obtained by cutting $\Gamma$ along $A$, so $s\geq 2$. Without loss of generality, we can assume that $Q \in \Gamma_1$. Let $\tilde \Gamma_1 \subset \Gamma_1$ be  the set of all the points $P$ of $\Gamma$ such that $P$ lies in the support of $D_P$. For the sake of a contradiction, suppose that $Q \notin \tilde\Gamma_1$.  Let $X$ be the closure  in $\Gamma$ of the connected component of $\Gamma_1 \setminus \tilde \Gamma_1$ which contains $Q$. Note that $X$ is a rational tropical curve and for all the points $P$ of  the boundary of $X$ in $\Gamma$, $P$ belongs to the support of $D_P$ (by the continuity of $\Red$). However, for all the points $P$ in the interior of $X$,   $P$ does not belong to the support of $D_P$. 
 We prove the following claim:
\begin{center}
\emph{Let $P$ be a point on the boundary of $X$. Then the $P$-reduced divisor $D_P$ is reduced for all the points of $X$. In particular all the points of the boundary of $X$ belong to the support of $D_P$ and no interior point lies in the support. 
}
\end{center}

This will certainly be a contradiction: for the point $Q$ in the interior of $X$, the cut $\Gamma \setminus \mathrm{int}(X)$ is saturated.  Indeed, all the points on the boundary of this cut have valence one, and the coefficient of $D_Q$ at these points is at least one since by the above claim all the points belong to the support of $D_Q$. 

\medskip 

To prove the claim, we use the explicit description of the reduced divisor map.  First of all, again as above, for any point $P$ on the boundary of $X$, there is a unique vector $\vec v$ emanating from $P$ in the direction of $X$. By the description of the reduced divisor in the interval $[P,P+\epsilon \vec v]$, $D_P$ is $(P+\epsilon\vec v)$-reduced for all $\epsilon$ small enough. Indeed, if this was not the case, then $P+\epsilon\vec v$ would belong to the support of $D_{P+\epsilon\vec v}$, which is certainly not true by the definition of $X$. 
For the sake of a contradiction, suppose now that the claim does not hold. Then there exists a point $R$ in the interior of $X$ and a direction $\vec v$ emanating from $R$ such that $D_R = D_P$ (for the boundary point $P$ of $X$ as above) such that $D_P$ is not reduced for any point $R+\epsilon \vec v$, for $\epsilon$ sufficiently small. In this case, $R$ should be in the support of $D_R$. Indeed, as we saw in the proof of Theorem~\ref{thm:continuity}, if the coefficient of $R$ in $D_R$ was zero, then  $D_R$ would remain obviously reduced for all the points sufficiently close to $R$, which is not the situation here, and so $D_R( R ) \geq 1$. But this is again a contradiction by the choice of $X$: no point $R$ in the interior of $X$ can belong the the support of  $D_R$.  Since $X$ is connected and $\Red$ is continuous, this shows that all the reduced divisors $D_R$, for $R \in  X$, are equal to $D_P$, and the claim follows.  This finishes the proof of the theorem.
 \end{proof}

\section{Maps into Tropical Projective Spaces}\label{sec:maps}
In this section, we give some complementary results on the image of the map $\Red$ and relate this map to maps into tropical projective spaces.

  Let $D$ be a divisor on a tropical curve $\Gamma$. Let $R(D)$ be the space of all rational functions $f$ such that $\div(f)+D \geq 0$. Since two rational functions $f$ and $g$ have the same divisor if and only if $f = c\odot g$ for a constant $c \in \mathbb T$, we have $|D| = \mathbb P R(D):= R(D) / \mathbb T$, the quotient is obtained by identifying $f$ with $c\odot f$. 

\subsection{The tropical module structure of $R(D)$} For two rational functions $f$ and $g$, the function $f\oplus g $ is a rational function on $\Gamma$ (recall that $(f\oplus g)(x) : = \max \{f(x),g(x)\}$). In addition, it is quite easy to show that $\div(f\oplus g) + D \geq 0$. It follows that $R(D)$ admits the structure of a tropical module over $\mathbb T$, the semi-field of tropical numbers. (Recall that the tropical multiplication $c\odot f$, for $c\in \mathbb T$ and $f \in R(D)$, is defined by 
$\bigl(c\odot f\bigr)(x) := f(x) +c$.)

Let $D$ be a divisor of non-negative rank. For a point $P$ in $\Gamma$, let $f_{P}$ be the unique rational function in $R(D)$ such that 
\begin{itemize}
\item[$(i)$] $D+\div(f_P) = D_P$, the unique $P$-reduced divisor equivalent to $D$.
\item[$(ii)$] $f_P(P) = 0$.
\end{itemize}

The fact that $f_P$ is unique comes from the unicity of $D_P$ and the assumption on the value of $f_P$ at $P$. We show below a set of properties of the rational functions $f_P$. Some complement properties are given in Section~\ref{sec:complement} below.

The first lemma gives a precise characterization of the rational function $f_P$. For a set of function $S=\{h:\Gamma \rightarrow \mathbb T\}$, the tropical sum $\bigoplus_{h \in S} h$ is the function defined by  $x \rightarrow \max_{h\in S}h(x)$ for all $x\in \Gamma$.

\begin{lemma}\label{lem:maxreduced2} Let $D$ be a divisor of non-negative rank on a tropical curve and $f_P$ be the rational function defined above. We have 
\[f_P = \bigoplus_{h \in R(D):\,\, h(P) = 0} h.\]

\end{lemma}

\begin{proof}
Since $f_P$ itself appears in the tropical sum in the right hand side of the above equation, we only need to show that for all $Q \in \Gamma$ and $h \in R(D)$ with $h(P) = 0$, we have $f_P(Q) \geq h(Q)$. By Lemma~\ref{lem:maxreduced}, we know that
 $h-f_P $ takes its maximum at $P$. This proves the lemma since this shows for all $Q \in \Gamma$, $h(Q) - f_P(Q)  \leq h(P)-f_P(P)=0$. 
\end{proof}

\begin{remark}\label{rem:module}\rm
As a tropical module, $R(D)$ cannot be generated in general  by the set $\bigl\{f_{P}\bigr\}$ for $P \in \Gamma$.
Indeed, for $g \in R(D)$ the existence of a solution $\{x_P\,\,; P\in \Gamma\}$ for the equation $g=\bigoplus_{P\in \Gamma} x_P\odot f_P$ (or, equivalently, for the family of equations $g(Q) = \max_{P\in \Gamma}\, x_P + f_P(Q)$ for all $Q$), can be  equivalently described as  follows. First note that the unique solution of the above equations, if exists, is given by 
$x_P = \min_{Q \in \Gamma} g(Q) - f_P(Q)$, and the following claim has  to be true:
\begin{center}
{\it if $S_P$ denotes the set of all points $Q$ such that $x_P = g(Q) -f_P(Q)$, then $\bigcup_{P\in \Gamma}S_P = \Gamma$.}
\end{center}
Let $E=\div(g)+D=\div(g-f_P)+\div(f_P)+D=\div(g-f_P)+D_P$. If $Q \in S_P$, and $g \neq c\odot f_P$ for all $P\in \Gamma$ and $c\in \mathbb T$, then $g-f_P$ takes at least two values, and it takes its minimum at $Q$, so for every point $v$ on the boundary of $S_P$, we have $E(v)\geq \deg^{out}_{S_P}(v)$. Call a cut $X$ good for $E$ if for all points $v$ on the boundary of $X$, we have  $E(v)\geq \deg^{out}_{X}(v)$. As we just noticed, the existence of a solution $\{x_P\,\,;P\in\Gamma\}$ implies that all the cuts $S_P$ are good for $E$ and they cover the entire $\Gamma$. However, there exist divisors which are not reduced for any point of $\Gamma$ and for which the union of the good cuts is not $\Gamma$. For example, consider  a tropical curve with two branching points of degree three and three parallel edges between them. Let $P$ and $Q$ be two points in the interior of two different edges. The divisor $E = 2P + 2Q$ is not reduced for any point of $\Gamma$ and the only cuts $X$ which are good for $E$ are $\{P\}$, $\{Q\}$, and $\{P,Q\}$. We do not know the answer to the following question:  {\it provided that the union of all the good cuts for $E$ is $\Gamma$, is it true that $g$ can be written as $g=\bigoplus_{P\in \Gamma} x_P\odot f_P$ for some $x_P$s  in  $\mathbb T\,?$}
\end{remark}

\subsection{On the image of the reduced divisor map} In this section we consider the image of the tropical curve $\Gamma$ under the map $\Red$ and study the basic properties of this map. We then relate our map $\Red$ to maps from $\Gamma$ to tropical projective spaces considered in~\cite{HMY09}. 
 
 We start by recalling the following basic result from~\cite{HMY09} (see the discussion in Section~\ref{sec:integeraffine})

\begin{proposition}[Hasse, Muskier, and Yu~\cite{HMY09}]\label{prop:HMY} Let $D'$ be an element of $|D|$ and $X$ be the set of all the points of the support of $D'$ which lie in the interior of an edge of $\Gamma$.  The dimension of the cell of $|D|$ which contains $D'$ is one less than the number of connected components of $\Gamma\setminus X$.
\end{proposition}

\begin{corollary}
The image of $\Red$ lies in the 1-skeleton of the cell complex $|D|$.
\end{corollary}
\begin{proof}
Let $E=D_P$ be the $P$-reduced divisor equivalent to $D$. Let $X$ be the set of all points in the support of $E$ which lie in the interior of an edge. Two case can happen:
\begin{itemize}
\item $P \notin X$. In this case we claim that $\Gamma \setminus X$ is connected, i.e., $E$ is a vertex of the cell complex $|D|$. Indeed, if this is not the case, then  taking a connected component of $\Gamma \setminus X$ which does not contain the point $P$, we obtain a cut in which all the boundary points are in the interior of an edge and lie in the support of $E$. This clearly contradicts the assumption that $E$ is $P$-reduced. 
\item $P \in X$. We claim that $\Gamma \setminus X$ contains either  one or two connected components. If this was not the case, then there would exist a connected component of $\Gamma \setminus X$ which does not contain $P$, and we again obtain a contradiction to the assumption that $E$ is $P$-reduced. 
\end{itemize}
This combined with Proposition~\ref{prop:HMY} above clearly completes the proof. 
\end{proof}

\begin{theorem}\label{thm:trop-reduced}
Let $\Gamma$ be a tropical curve. Every point in the image of $\Red$ is either a vertex of $|D|$ or a tropical convex combination of two vertices of $|D|$ which are in the image of the map $\Red$.
The tropical convex hull of the image of $\Red$ is finitely generated. \end{theorem}
\begin{proof}
Let $T$ be the set of all the points $D_Q$ which are among the vertices of $|D|$. We first claim that the tropical convex hull of the image of $\Red$ is generated by the set of all $D_Q \in T$. For this aim, we will only need to show that for any point $P\in \Gamma$, the $P$-reduced divisor $D_P\sim D$ can be written as a convex hull of two elements of $T$. Let $P\in  \Gamma$. If $D_P \in T$, we do not have anything to prove. So suppose that $D_P \notin T$. This means that the set $X$ consisting of all the points in the support of $D_P$ which lie in the interior of an edge contains the point $P$ and $\Gamma \setminus X$ consists of two connected components.  Let $\vec u$ and $\vec w$ be the two unit vectors tangent to $\Gamma$ at $P$.  By the explicit description of the reduced divisor map we gave in Section~\ref{sec:reducedmap}, it is easy to see that $f_P $ can be written as the tropical convex hull of the rational functions $f_{P+\epsilon_0 \vec u}$ and $f_{P+\delta_0\vec w}$ where $\epsilon_0$ and $\delta_0$ are respectively the supremum value of $\epsilon$ and $\delta$ such that for all the values of $\epsilon'$ in $[0,\epsilon]$ and $\delta'$ in $[0,\delta]$,  the support points of the corresponding reduced devisors living in the interior of the edges form a cut. The two divisors $D_{P+\epsilon_0 \vec u}$ and $D_{P+\delta_0\vec w}$ are certainly vertices of $|D|$. We conclude that every point in the image of $\Red$ is a tropical convex combination of two vertices of $|D|$ which lie in the image of $\Red$. This also shows that the tropical convex hull of the image of $\Red$ is a finitely generated module. 
\end{proof}

We define the reduced linear system of $D$ as follows. 

\begin{definition}\rm
Define $\tilde R(D)$ to be the tropical convex hull of the image of $\Red$, i.e., the tropical module generated by the rational functions $f_P$. The {\it reduced linear system} of $D$, denoted by $|D|_r$ is defined to be $\tilde R(D)/\sim$.
\end{definition}

Contrary to what one might expect, it is not  in general true that the tropical convex hull of the image of $\Red$ has dimension equal to $r(D)$. An example is given by a tropical curve $\Gamma$ consisting of two branching points of degree three $P$ and $Q$ connected with three parallel edges. The divisor $2P+Q$ has rank one while the tropical convex hull of the image of $\Red$ has dimension two. 

\subsection{Maps to tropical projective spaces} \label{sec:map-projective}
Given a fixed model for $\Gamma$, by Theorem~\ref{thm:trop-reduced}, $\tilde R(D)$ is generated by all the vertices of $|D|$ which are reduced with respect to some base point in $\Gamma$. Let $\mathcal F$ be the set of all the rational functions $f_P$ corresponding to these vertices. $\mathcal F$ allows to define a map from $\Gamma$ to $\mathbb T\mathbb P(\mathcal F^*)$. This later space is the tropical projective space consisting of all the maps $\mathcal F \rightarrow \mathbb T$ modulo tropical multiplication by a constant. Note that the map $\phi_{\mathcal F}:\Gamma \rightarrow \mathbb T\mathbb P(\mathcal F^*)$ is canonical once a model for $\Gamma$ is fixed. In particular, if $\Gamma$ is of genus at least two and the model of $\Gamma$ consists of the coarsest model (i.e., all the vertices are the branching points), then the map $\Gamma \rightarrow \mathbb T\mathbb P(\mathcal F^*)$ is canonical.

The following theorem establishes a precise relation  between our map $\Red$ and the embeddings of tropical curves into tropical projective spaces considered in~\cite{HMY09}.

\begin{theorem}\label{thm:embedding2}
\begin{itemize}
\item[1.] The image of $\phi_{\mathcal F}$ is homeomorphic to the image of $\Red$.
\item[2.] In addition, there is an injective map $|D|_r \rightarrow \mathbb{TP}(\mathcal F^*)$ with image $\textrm{trop-conv}\Bigl(\phi_{\mathcal F}(\Gamma)\Bigr)$, the tropical convex hull of the image of $\phi_{\mathcal F}$ in $\mathbb T\mathbb P({\mathcal F}^*)$. In particular, $|D|_r$ is homeomorphic to $\textrm{trop-conv }\Bigl(\phi_{\mathcal F}(\Gamma)\Bigr)$. In other words, there exists  a commutative diagram

\end{itemize}
\[
   \begin{diagram}
     \node{\Gamma} \arrow{e,t}{\Red} \arrow{s,r}{id} 
       \node{|D|_r}  \arrow{s,r}{\sim} \\
    \node{\Gamma} \arrow{e,t}{\phi_{\mathcal F}} \node{\textrm{trop-conv }\Bigl(\phi_{\mathcal F}(\Gamma)\Bigr)} \arrow{e,t,J}{\iota}\node{\mathbb T \mathbb P(\mathcal F^*)}
   \end{diagram}
\]
\end{theorem}
\begin{proof}
\begin{itemize}
\item[1.] Define a map $\Psi: \phi_{\mathcal F}(\Gamma) \rightarrow |D|_r = \tilde R(D)/\sim $ as follows. Given a point $P \in \Gamma$, let $\Psi\bigl(\phi_{\mathcal F}(P)\bigr) = \oplus_{f \in \mathcal F} (-f(P))\odot f$. We claim that $\Psi\bigl(\phi_{\mathcal F}(P)\bigr) = f_P$ in $\tilde R(D)$. From the description we give below, this will prove the first part of the theorem. To show this, note that since $\mathcal F$ generates $\tilde R(D)$ as a tropical module, there are $x_f \in \mathbb T$ for $f \in \mathcal F$ such that $\oplus_{f\in \mathcal F}\, x_f \odot f =f _P$. This means that $x_f+f(P) \leq f_P(P)=0$ for all $f \in \mathcal F$, or equivalently $x_f \leq -f(P)$.  In particular, we must have $f_P=\oplus_{f\in \mathcal F}\, x_f \odot f \leq \Psi\bigl(\phi_{\mathcal F}(P)\bigr)$. Since $(-f(P))\odot f$ takes value zero at $P$, by applying Lemma~\ref{lem:maxreduced2}, we have $f_P = \Psi\bigl(\phi_{\mathcal F}(P)\bigr)$ and the claim follows. To prove the first part, note that there is a map from $\Red(\Gamma) \rightarrow \phi_{\mathcal F}(\Gamma)$ obtained by sending $f_P$ to $\phi_{\mathcal F}(P)$. This map is well-defined since $f_P=\oplus_{f \in \mathcal F}\,\, x_f\odot f $ and $-f(P)$ is the maximum value $x_f$ can take while verifying the above equation. Thus, for two points $P$ and $Q$ with $f_P =f_Q$, we have $f(P) =f(Q)$ for all $f \in \mathcal F$.

\item[2.]  Define a map $\Theta: |D|_r \rightarrow \mathbb T \mathbb P(\mathcal F^*)$ as follows. Given an element $g \in |D|_r$, $\Theta(g) \in \mathbb T \mathbb P(\mathcal F^*)$ is the map which sends $f \in \mathcal F$ to $\max_{Q \in \Gamma} f(Q)-g(Q)$. 
\[\Theta(g): \mathcal F \rightarrow \mathbb T\,\,\textrm{ is given by } \Theta(g)(f) = \max_{Q\in \Gamma} f(Q) - g(Q).\]
As mentioned in Remark~\ref{rem:module}, an element $g \in |D|_r = \textit{trop-conv}\Bigr(\{f_P\}_{P \in \Gamma}\Bigr)$ can be written in the form $g=\oplus_{f \in F}\, x_f \odot f$ where $x_f := \min_{Q \in \Gamma} g(Q)-f(Q)$, so that the map $\Theta(g)$ sends $f$ to $-x_f$. By definition, it is clear that $\Theta(\lambda\odot g) = (-\lambda)\odot \Theta(f)$, so that $\Theta$ is a well-defined map from $|D|_r$ to $\mathbb T\mathbb P(\mathcal F^*)$. Note also that the restriction of $\Theta$ on $\Red(\Gamma)$ is the inverse of the map $\Psi$ constructed above. Indeed if $g =f_P$ for a point $P \in \Gamma$, then  the rational function $f-f_P$ takes its maximum value at $P$, and this value is equal to $f(P)$. Thus, $\Theta(f_P)=\phi_{\mathcal F}(P).$ In addition, the image of $\Theta$ clearly lies in the tropical convex hull of the image of $\phi_{\mathcal F}$. By the very definition, the injectivity of $\Theta$ follows. To show that $\Theta$ is surjective, let $\eta: \mathcal F \rightarrow \mathbb R$ be  a given map sending $f$ to $\eta_f$ and $\eta = \oplus_{Q\in \Gamma}\, h(Q)\odot\phi_{\mathcal F}(Q)$. Define the rational function $g$ in $|D|_r$ as the tropical sum $g:= \oplus_{f \in \mathcal F}\, (-\eta_f) \odot f$. We claim that $\Theta(g) =\eta$. 

By assumption, we have $\eta_f = \max_Q\, h(Q)+f(Q)$. Thus, $-h(Q) \geq f(Q) - \eta_f$ for all $f$, i.e., 
\begin{align*}
-h(Q) \geq g(Q),\qquad \forall \,\, Q\in\Gamma. 
\end{align*}
 Note also that by the definition of $g$, we have for fixed $f$, 
\begin{align*}
\eta_f \geq f(Q) - g(Q) \qquad \forall\,\, Q\in\Gamma.
\end{align*}
Since $\eta_f =\max_{Q\in \Gamma}\, h(Q)+f(Q)$, there is a $Q_0$ such that $\eta_f = h(Q_0)+f(Q_0)$. For this particular point $Q_0$, the above inequality becomes 
\[
h(Q_0)+f(Q_0) \geq f(Q_0) -g(Q_0), \qquad \mathrm{i.e.,}\qquad -h(Q_0) \leq g(Q_0). 
\]

We infer that $h(Q_0) = -g(Q_0)$, and $\eta_f = f(Q_0) -g(Q_0)$.

This shows that for all $f$, there is a point $Q$ such that $-h(Q) = g(Q) = f(Q) -\eta_f$, or $\eta_f = f(Q)-g(Q)=\max_Q\,f(Q)-g(Q) =\Theta(g)(f)$.  \end{itemize}
\end{proof}
\begin{corollary}
Let $\mathcal F$ be defined as above. Then $D$ is very ample if and only if $\phi_{\mathcal F}$ is injective. 

\end{corollary}

\subsection{Duality}\label{sec:complement}
Let $f_P$ be as before the rational function defining the $P$-reduced divisor $D_P$, i.e.,  $D_P=D + \div(f_P)$, and $f_P(P)=0$. Consider the function $f_Q^t$ obtained from the family of rational functions $f_P$ by interchanging the role of base points and arguments, namely, $f_Q^t$ is the function defined by $f_Q^t(P) = - f_P(Q)$. We now show that all the functions $f^t_Q$ are rational and their tropical convex hull contains the tropical module $R(D)$. (Note that $f^t_Q$ does not necessarily lie in $R(D)$.)

To show the second assertion, note that for a point $P \in \Gamma$, we have $g(P)\odot f_P^t(P) = g(P)+0 = g(P)$, so we  only need to show that for all $P$ and $Q$ in $\Gamma$, $g(P)\odot f_P^t(Q) \leq g(Q)$, or, $g(P)- f_Q(P) \leq g(Q) - f_Q(Q).$ This directly follows from the fact that $g-f_Q$ takes its maximum value at $Q$, c.f.~Lemma~\ref{lem:maxreduced}.

The fact that $f^t_Q$ is a rational function is an easy consequence of the results of the previous sections.  Below, we give a more explicit proof which provides a description of the sectional derivatives of $f^t_Q$ in terms of the results of Section~\ref{sec:reduceddivisormap} .  For fixed $Q\in \Gamma$, and for any $P\in \Gamma$ and any unit vector $\vec u$ emanating from $P$,  we need to calculate $\frac d{d\epsilon}\, f^t_Q(P+\epsilon \vec u)\,|_{\epsilon=0}$, and show that it is an integer. By definition, this is $- \frac d{d\epsilon}\, f_{P+\epsilon \vec u}(Q)\,|_{\epsilon=0} = \lim_{\epsilon \rightarrow 0} \frac {- f_{P+\epsilon \vec u}(Q) + f_{P}(Q)}{\epsilon}$. Note that $g_\epsilon(.) = f_{P+\epsilon \vec u}(.) - f_{P}(.)$ satisfies $D_P +\div(g_\epsilon) = D_{P+\epsilon \vec u}$. This shows that $g_\epsilon$ is the rational function given in the proof of Theorem~\ref{thm:continuity} up to an additive constant. The constant is defined by the requirement that $f_{P+\epsilon \vec u}(P+\epsilon \vec u)=0$. A simple calculation shows that the derivation is exactly
\[- \frac d{d\epsilon}\, f_{P+\epsilon \vec u}(Q)\,|_{\epsilon=0} = \frac d{d\epsilon}\,f_P(P+\epsilon \vec u) - h_{\vec u}(Q), \]
where the function $h_{\vec u}$ is defined as follows. Let $Y$ be the cut considered in the proof of Theorem~\ref{thm:continuity} which contains $P$, does not contain $P+\epsilon \vec u$ and which is maximum with respect to the property that for all $v$ on its boundary, $D_{P}(v) \geq \deg^{out}_Y(v)$. If such a cut does not exist simply let $Y = \emptyset$.  Then $h_{\vec u}(Q) = 0$ if $Q$ lies outside $Y$ and $h_{\vec u}(Q) = D_{P}(P) - \deg^{out}_Y(P)+1$ if $Q$ is in  $Y$. In particular, the slope of $f^t_Q$ at $P$ is an integer. It is clear that the number of these slopes is finite. This shows that $f^t_Q$ is a rational function.  By the description of the directional derivatives of $f^t_Q$ at $P$, we also see that the coefficient of $D+ \div(f^t_Q)$ at $P$ is $D_{P}( P ) - \sum_{\vec u} h_{\vec u} (Q)$, where the sum is over all the unit vectors $\vec u$ tangent to $\Gamma$ at $P$, and $h_{\vec u}(Q)$ is either zero or equal to $D_{P}( P ) - \deg^{out}_Y(P)+1$ depending on whether $Q$ lies outside or inside $Y$, respectively.

\section{Algebraic Curves and the Reduced Divisor Map}\label{app:weier}
To make the analogy between the classical and the tropical case,  we provide in this section the classical analogue of the reduced-divisor map.

\medskip 

Let $C$ be a non-singular projective curve over an algebraically closed field $\kappa$ of characteristic zero. Let  $D$ be a divisor on $C$ and $P$ be a point of $C$. A divisor $D'\sim D$ is called $P$-reduced if all the coefficients of $D'$ outside $P$ are non-negative, and in addition, the coefficient of $P$ in $D'$ is maximum among all the divisors linearly equivalent to $D$ with this property.   
\begin{proposition}
For every divisor $D$ on $C$, there exists a unique $P$-reduced divisor linearly equivalent  to $D$.
\end{proposition}
\begin{proof} The existence of divisors $D'\sim D$ with non-negative coefficients outside $P$ is easy (e.g., is a direct consequence of Riemann-Roch). Since the coefficients of $D'$ are bounded above by the degree of $D$, the existence of $P$-reduced divisors follows. To show the uniqueness, note that if $f$ and $g$ are two non-collinear rational functions such that both $D' = D +\div(f)$ and $D'' = D+\div(g)$ are $P$-reduced, and so  have the same coefficient at $P$, then for two appropriately chosen non-zero  scalars $\alpha$ and $\beta$ in $\kappa$,  the non-zero rational function $\alpha f + \beta g$ has an order at $P$ strictly larger than that of $f$ and $g$. In addition $D + \div(\alpha f + \beta g)$ has non-negative coefficients outside $P$. This leads to a contradiction and the proposition follows. 
\end{proof}

 Let $D$ be a divisor on $C$ of non-negative rank, and let $\mathcal L = \mathcal L(D)$ be the corresponding invertible sheaf on $C$. Let $L$ be the vector space of the global sections of $\mathcal L(D)$; $L$ is the vector space of all the rational functions $f$ with $D+\div(f) \geq 0$ ($L=R(D)$ with our earlier notation). Note that $\dim (L) = r(D)+1$. Denote by $L^*$ the dual of $L$.

   \vspace{.3cm}
   
 There are two maps from $C$ to a projective space defined by the above data that we now explain. 
 
 \noindent The easier to define, which is the one usually considered, is the map $\phi: C \rightarrow \mathbb P(L^*)$: $\phi(P)$ is the point of $\mathbb P(L^*)$ corresponding to the hyperplane of $L$ defined by all the rational functions $f \in L$ such that the coefficient of $P$ in $D+\div(f)$ is at least one. Equivalently, the map $\phi$ is the map given generically by evaluating at the points of $C$, i.e., generically, $\phi(P)$ is the linear form on $L$ defined by sending a rational function $f \in L$ to $f(P)$. If $D$ is base-point free, $\phi$ is well-defined globally. In any case, the map can be extended to whole of $C$.
  
  \vspace{.2cm}
  
 \noindent The second map is $\eta : C \rightarrow \mathbb P(L)$ defined by reduced divisors. Suppose that $D$ has non-negative rank. For every point $P \in C$, let $D_P$ be the unique $P$-reduced divisor linearly equivalent to $D$. Since $r(D) \geq 0$, we have $D_P \geq 0$. Let $f_P$ be a rational function in $L$ such that $D_P = D+\div(f_P)$. Note that modulo multiplication by a constant,  $f_P$ is well-defined. In other words, the line $\ell_P$ defined by  $f_P$ is well-defined in $L$.  Define a map $\eta : C \rightarrow \mathbb P(L)$ by setting $\eta(P) = \ell_P$ for any point $P \in C$. 

\noindent We will show below that $\eta$ is a morphism defined by a set of sections of $\mathcal L^{\otimes r}\,\Omega^{\otimes r(r-1)/2}$, where $\Omega$ is the invertible sheaf of differentials on $C$ and $r=r(D)$. These sections are defined in terms of Wronskians.

 \subsection{Wronskians and the reduced divisor map}
 
 We first  obtain a generic description of the map $\eta$. Let $\mathbb K=\mathbb K( C )$ be the field of rational functions on $C$. 
 Choose a basis $f_0,\dots,f_r$ of $L$ over $\kappa$ and let $\mathcal F = \{f_0,\dots,f_r\}$. Let $\Omega$ be the invertible sheaf of differentials  on $C$. Let $P$ be a point of $C$ and $\O_P$ be the local ring of the structural sheaf $\O_C$ at $P$. $\Omega_P$ is the module of $\kappa$-differentials of $\O_P$ and $d:\mathcal O_O \rightarrow \Omega_P$ is the derivation.   
 
\noindent  If $\tau$ denotes a local parameter at $P$, a basis for $\Omega_P$ is given by $d\tau$. 
 
\noindent   The $\O_P$-module $\mathcal L_P$ is free of rank one, suppose that it is generated by $g_P\in \mathcal L_P$.   

\medskip 

\noindent  For each $i \in \{0,\dots, r\}$, there exists an element $f_{i,P}\in \O_P$ such that $f_{i|P} = f_{i,P}g_P$. 
 
 For each $i \in \{0,\dots,r\}$, inductively define $f^{(j)}_{i,P}\in \mathcal O_P$ as follows. Let $f^{(0)}_{i,P} := f_{i,P}$, and suppose that $f^{(j)}_{i,P}\in \O_P$ is already defined. Then $f^{(j+1)}_{i,P}:= df^{(j)}_{i,P}/d\tau \in \O_P$ (i.e., $f^{(j+1)}_{i,P}$ is the function defined by $df^{(j)}_{i,P} = f^{(j+1)}_{i,P} d\tau$).
 
 \medskip 
 
 The Wronskian of $f_0,\dots,f_r$ at $P$, denoted by $\Wr_{\mathcal F,P}$, is now defined as 
 $$\Wr_{\mathcal F,P} = \det\Bigl(f^{(j)}_{i,P}\Bigr)_{i,j =0}^{r} g_P^{r+1} \bigl(d\tau\bigr)^{r(r+1)/2} \in \Gamma \bigl(\O_P,\mathcal L_P^{\otimes(r+1)} \Omega^{\otimes r(r+1)/2}_{P}\bigr).$$
 It is easy to show that the above (local) definition does not depend on the choice of the local parameter $\tau$ and glue together to define a non-zero global section $\Wr_{\mathcal F}$ of $\mathcal L^{\otimes(r+1)} \Omega^{\otimes r(r+1)/2}$, i.e., $\Wr_{\mathcal F} \in \Gamma\bigl(C, \mathcal L^{\otimes(r+1)} \Omega^{\otimes r(r+1)/2}\bigr)$. The fact that $\Wr_{\mathcal F}$ is non-zero comes from the non-triviality of the determinant $\det\Bigl(f^{(j)}_{i,P}\Bigr)_{i,j =0}^{r}$ for linearly independent $f_{i,P}$'s. This later claim can be seen easily by looking at the completion $\hat \O_P\sim \kappa[[\tau]]$ of $\O_P$ at the maximal ideal $m_P$ of $P$ and the map $\O_P \rightarrow \kappa[[\tau]]$.   The image of $\Wr_{\mathcal F,P}$ in $\kappa[[\tau]]$ is non-zero, see for example~\cite{BD10} for a simple proof. 
 
 Note that in the discussion above we only used the linear independence of the sections $f_i$, thus, it is clear that instead of choosing a basis for $L$, we can take any family $\mathcal G$ consisting of $s \leq r+1$ linearly independent elements of $L$ and define the Wronskian $\Wr_\mathcal G$ with respect to that family. The Wronskian $\Wr_{\mathcal G}$ will be then a non-zero global section of $\mathcal L^{\otimes s} \omega^{\otimes s(s-1)/2}$. 
 
 \medskip

 For $i\in \mathcal F$, define $\mathcal F_i = \mathcal F \setminus \{i\}$ and consider the Wronskian $\Wr_{\mathcal F_i} \in \Gamma\bigl(C, \mathcal L^{\otimes r} \Omega^{\otimes r(r-1)/2}\bigr)$. Locally at a point $P$, and in the terminology of the previous paragraph, $\Wr_{\mathcal F_i}$ is the maximal minor of the matrix $\Bigl(f^{(j)}_{i,P}\Bigr)_{i=0,\dots, r, j =0,\dots, r-1}$ obtained by deleting the $i^{\,\textrm{th}}$ column. The sections $\Wr_{F_i}$ satisfy the following equations: 
 
 \begin{align}\label{eq:wronsk}
 \forall \,\, 0\leq j \leq r-1 \qquad &\sum_{i=0}^r (-1)^i\Wr_{\mathcal F_i,P}\, f^{(j)}_{i,P} = 0, \qquad \textrm {and}\\
  &\label{eq:wronsk2} \sum_{i=0}^r (-1)^i\Wr_{\mathcal F_i,P}\, f^{(r)}_{i,P} = \Wr_{\mathcal F,P}\, .
 \end{align}

Since $\Wr_{\mathcal F} \neq 0$, by the definition of $f^{(j)}_{i,P}$,  it follows that generically the rational function $\sum_{i} (-1)^i\Wr_{\mathcal F_i,P}(P) f_i$ is exactly the one defining the reduced divisor $D_P$, i.e., the line defined by $\sum_{i} \Wr_{\mathcal F_i,P}(P) f_i$ is $\ell_P$ (the line defined by the rational function $f_P$ with $D+\div(f_P)$ $P$-reduced).   Here $\Wr_{\mathcal F_i,P}(P)$ is the value of $\Wr_{\mathcal F_i,P}$ at $P$. 

\vspace{.2cm}

\noindent Indeed, to see this, write $f_P = \sum a_i f_i$ for $a_i \in \kappa$, and note that $f_P$ has a zero of order at least $s$ at $P$ if and only if all the derivations of order at most $s-1$ of $f_P$ at $P$ are zero, i.e., $\sum_i a_i f_{i,P}^{(j)}(P) =0$. Now use $\Wr_F \neq 0$ and Equations~(\ref{eq:wronsk}) and~(\ref{eq:wronsk2}) to obtain the claim.
 
 \vspace{.2cm} We infer that, generically, the map $\eta$ coincides with the map to $\mathbb P(L) \, \Bigl(\sim \mathbb P^r$ with respect to the basis $f_0,\dots,f_r\Bigr)$ defined by the sections $\{(-1)^i\Wr_{\mathcal F_i}\}_{i=0}^r$ of $\mathcal L^{\otimes r}\Omega^{\otimes r(r-1)/2}$ $\Bigl($i.e., $\eta$ sends the point $P$ to the point $[(-1)^0\Wr_{\mathcal F_0}(P):\dots:(-1)^r\Wr_{\mathcal F_r}(P)]$ of $\mathbb P^r \Bigr)$. Let $\overline \Wr$ be the extension of the map defined by the sections $\Wr_{\mathcal F_i}$ to the whole curve $C$.
 
 \medskip 
 We now show
  \begin{theorem}\label{thm:app-clas}
 The map $\eta$ defined by the reduced divisors coincides globally with the morphism $\overline \Wr$ defined by the sections $(-1)^i\Wr_{\mathcal F_i}$ of $\mathcal L^{\otimes r}\Omega^{\otimes r(r-1)/2}$, $i=0,\dots,r$. 
  \end{theorem}
 \begin{proof}
 We have already seen that $\eta$ generically coincides with $\overline \Wr$, e.g., for those points of $C$ for which the determinants appearing in the definition of the Wronskians $\Wr_{\mathcal F_i}$ are not all zero. For a point $P$ and $i \in \{0,\dots,r\}$, let $A_{i,P}$ be the determinant in the definition of $\Wr_{\mathcal F_i}$, in the terminology of the previous lines:
 $$A_{i,P} = \det\Bigl(f^{(j')}_{j,P}\Bigr)_{j=0,\dots,\hat i,\dots, r}^{j' = 0,\dots, r-1}.$$ 
 We now show that the image by $\overline \Wr$ of a point $P$ for which all the determinants $A_{i,P}$ are zero  is exactly $\eta(P)$. This will finish the proof of the theorem.

 Let $P$ be such a point. And let $v\geq 1$ be the integer such that $A_{i,P} = \tau^{v} h_i$ for all $i$, and at least for one $i$, $h_i(P) \neq 0$. The image of $P$ by $\overline \Wr$ is the point $[h_0(P):-h_1(P):\dots:(-1)^rh_r(P)]$ of $\mathbb P^{r} \Bigl(= \mathbb P(L)$ in the basis $\mathcal F\Bigr)$. We have to show that $\sum (-1)^i h_i(P) f_{i,P}$ has the maximum order of zero at $P$. 
 
 \noindent By considering the injective map $\mathcal O_P \rightarrow \kappa[[\tau]]$, we can suppose that all the elements $f_{i,P}$ and $h_{i}$ are in $\kappa[[\tau]]$, and evaluating at $P$ consists in setting $\tau =0$.

 \noindent Write $f_{i,P} = \sum_{s \geq 0} a_{i,s} \tau^j$.  And consider the vectors $V_s = \bigl(a_{0,s},\dots,a_{r,s}\bigr)$. For $i=1,\dots,r$, let $s_i$ be the smallest index $s$ such that the span of the vectors $V_0,\dots,V_{s_i}$ has rank $i$. 
 In particular, the span of $V_{s_1},\dots, V_{s_r}$ has rank $r$. We now show that $h_0(P),\dots,h_r(P)$ are a constant multiple of the maximal minors of the $r(r+1)$ matrix $M$ with lines equal to the vectors $V_{s_1},\dots,V_{s_r}$,
 
 $$
M:= \left(
\begin{array}{cc}
V_{s_1} \\
 \vdots\\
 V_{s_r}
 \end{array}
\right)
 = \left(
\begin{array}{cccc}
a_{s_1,0} & a_{s_1,1} & \cdots & a_{s_1,r} \\
 \vdots & \vdots & &\vdots\\
a_{s_r,0}&a_{s_r,1}& \cdots& a_{s_r,r}
 \end{array}
\right)
$$

  This will of course finish the proof. Indeed,  if $\mu_i$ denotes the maximal minor obtained by deleting the $i^{\textrm{th}}$ column, then by the definition of $s_r$,  at least one of the $\mu_i$'s is non-zero, and in addition for any vector $V_s$
 
 $$\sum_{j=0}^r (-1)^i \mu_i a_{i,s} = 0\qquad \textrm{if and only if }\,  V_s \textrm{ belongs to the span of } V_{s_1},\dots,V_{s_r},$$
 thus,  $\sum_i (-1)^i\mu_i f_{i,P}$ has the maximum order of zero at $P$. 

 \vspace{.2cm}
  
  \noindent   We claim that $v = s_1+(s_2-1)+\dots+(s_r-r)$ and the coefficient of $t^v$ in $A_{i,P}$ is $(\prod_{i=1}^r \frac {s_i!}{i!})\,\, \mu_i$ (this is trivially equivalent to  $h_i(0) =(\prod_{i=1}^r \frac {s_i!}{i!}) \,\, \mu_i$).   The claim is easily obtained by developing the determinant in the definition of $A_{i,P}$. 
 For any $s$, let $V_s^i = \bigl(a_{s,0},\dots,\hat a_{s,i}, \dots, a_{s,r}\bigr)$, so that
 $$
\mu_i= \det\left(
\begin{array}{cc}
V_{s_1}^i \\
 \vdots\\
 V_{s_r}^i
 \end{array}
\right).
 $$
 
 We have 
 \begin{align*}
 A_{i,P}& = \det \Bigl(f_{j,P}^{(j')}\Bigr)_{j=0,\dots,\hat i,\dots, r}^{ j' = 0,\dots, r-1}
 = \det\left(
\begin{array}{cc}
\sum_{j_1\geq 0} \tau^{j_1} V_{j_1}^i \\
 \\
 \vdots\\
 \\
 \sum_{j_i\geq i-1} \frac{j_i !}{i!} \tau^{j_i-i+1}V_{j_2}^i\\
 \\
\vdots\\
\\
 \sum_{j_r \geq r-1} \frac {j_r !}{r!} \tau^{j_r-r+1} V_{j_r}^i
 \end{array}
\right)\\
 &
 \\
 &= \sum_{j_1\geq 0,\, j_2\geq 1,\,\dots,\, j_{r} \geq r-1} \,\, \tau^{j_1+(j_2-1)+\dots+(j_r-r+1)} \,\, \prod_{i=1}^r \frac {j_i!}{i!}\,\,\det\left(
\begin{array}{cc}
V_{j_1}^i \\
 \vdots\\
 V_{j_r}^i
 \end{array}
\right).
 \end{align*}
 By the definition of $s_1,\dots,s_r$, all the determinants 
 $
 \det \left(\begin{array}{cc}
V_{j_1}^i \\
 \vdots\\
 V_{j_r}^i
 \end{array}
\right)
 $
 are zero for $\sum_{i=1}^r(j_i - i+1) < \sum_{i=1}^r(s_i - i+1)$, and the coefficient of $\tau^{\sum_{i=1}^r(s_i - i+1)}$ is 
 exactly  $\prod_{i=1}^r \frac {s_i!}{i!}\,\, \mu_i$. The proof of the theorem is now complete.  
   \end{proof}

\medskip
 
Note that $D$-Weierstrass points are precisely the points  in the support of the effective divisor defined by the global section $\Wr_{\mathcal F}$ of $\mathcal L^{\otimes (r+1)}\Omega^{\otimes r(r+1)/2}$, and the multiplicity of a Weierstrass point is  defined to be the coefficient of  the point in this divisor. This is a simple consequence of the definition of $\Wr_{\mathcal F}$ and Equations~\eqref{eq:wronsk} and~\eqref{eq:wronsk2} above. The total number of $D$-Weierstrass points, counted with multiplicity, is given by the degree of the invertible sheaf $\mathcal L^{\otimes (r+1)}\Omega^{\otimes r(r+1)/2}$ which is $(r+1)\deg(D) +  r(r+1)(g-1)$.

\section*{Acknowledgment}
 Part of this work was done during a visit of the author at Charles University in Prague in 2008. The author is very grateful to Dan Kr\`al' and the members of KAM for the warm hospitality and great discussions.  Special thanks to Matt Baker for his interest in this work, for bringing the author's attention to the work of~\cite{HMY09} related to the results of this paper, and to the work of~\cite{Luo09}, and for his many constructive suggestions and interesting comments  on the first draft of this paper which helped to improve the presentation.   
 Many thanks to Erwan Brugall\'e for several interesting and fruitful discussions, and for motivating the author to write up these results. Special thanks to Pooran Memari and Mathieu Desbrun for their help and support during the preparation of this paper.

\bibliographystyle{amsplain}
\bibliography{transams}

\end{document}